
%
%
%
\def\unredoffs{} \def\redoffs{\voffset=-.31truein\hoffset=-.59truein}
\def\speclscape{\special{ps: landscape}}
%
%
%
%
\newbox\leftpage \newdimen\fullhsize \newdimen\hstitle \newdimen\hsbody
\tolerance=1000\hfuzz=2pt
\catcode`\@=11 
%
\ifx\answ\bigans\message{(This will come out unreduced.}
\magnification=1200\unredoffs\baselineskip=16pt plus 2pt minus 1pt
\hsbody=\hsize \hstitle=\hsize 
\else\message{(This will be reduced.} \let\l@r=L
\magnification=1000\baselineskip=16pt plus 2pt minus 1pt \vsize=7truein
\redoffs \hstitle=8truein\hsbody=4.75truein\fullhsize=10truein\hsize=\hsbody
\output={\ifnum\pageno=0 
  \shipout\vbox{\speclscape{\hsize\fullhsize\makeheadline}
    \hbox to \fullhsize{\hfill\pagebody\hfill}}\advancepageno
  \else
  \almostshipout{\leftline{\vbox{\pagebody\makefootline}}}\advancepageno
  \fi}
\def\almostshipout#1{\if L\l@r \count1=1 \message{[\the\count0.\the\count1]}
      \global\setbox\leftpage=#1 \global\let\l@r=R
 \else \count1=2
  \shipout\vbox{\speclscape{\hsize\fullhsize\makeheadline}
      \hbox to\fullhsize{\box\leftpage\hfil#1}}  \global\let\l@r=L\fi}
\fi
%
\newcount\yearltd\yearltd=\year\advance\yearltd by -1900

\def\Title#1#2{\nopagenumbers\abstractfont\hsize=\hstitle\rightline{#1}%
\vskip 1in\centerline{\titlefont #2}\abstractfont\vskip .5in\pageno=0}
%

%
%

\def\draftmode{\message{ DRAFTMODE }\def\draftdate{{\rm preliminary draft:
\number\month/\number\day/\number\yearltd\ \ \hourmin}}%
\headline={\hfil\draftdate}\writelabels\baselineskip=20pt plus 2pt minus 2pt
 {\count255=\time\divide\count255 by 60 \xdef\hourmin{\number\count255}
  \multiply\count255 by-60\advance\count255 by\time
  \xdef\hourmin{\hourmin:\ifnum\count255<10 0\fi\the\count255}}}
\def\nolabels{\def\wrlabeL##1{}\def\eqlabeL##1{}\def\reflabeL##1{}}
\def\writelabels{\def\wrlabeL##1{\leavevmode\vadjust{\rlap{\smash%
{\line{{\escapechar=` \hfill\rlap{\sevenrm\hskip.03in\string##1}}}}}}}%
\def\eqlabeL##1{{\escapechar-1\rlap{\sevenrm\hskip.05in\string##1}}}%
\def\reflabeL##1{\noexpand\llap{\noexpand\sevenrm\string\string\string##1}}}
\nolabels
%
\global\newcount\secno \global\secno=0
\global\newcount\meqno \global\meqno=1
\def\newsec#1{\global\advance\secno by1\message{(\the\secno. #1)}
\global\subsecno=0\eqnres@t\noindent{\bf\the\secno. #1}
\writetoca{{\secsym} {#1}}\par\nobreak\medskip\nobreak}
\def\eqnres@t{\xdef\secsym{\the\secno.}\global\meqno=1\bigbreak\bigskip}
\def\sequentialequations{\def\eqnres@t{\bigbreak}}\xdef\secsym{}
\global\newcount\subsecno \global\subsecno=0
\def\subsec#1{\global\advance\subsecno by1\message{(\secsym\the\subsecno. #1)}
\ifnum\lastpenalty>9000\else\bigbreak\fi
\noindent{\it\secsym\the\subsecno. #1}\writetoca{\string\quad
{\secsym\the\subsecno.} {#1}}\par\nobreak\medskip\nobreak}
\def\appendix#1#2{\global\meqno=1\global\subsecno=0\xdef\secsym{\hbox{#1.}}
\bigbreak\bigskip\noindent{\bf Appendix #1. #2}\message{(#1. #2)}
\writetoca{Appendix {#1.} {#2}}\par\nobreak\medskip\nobreak}
%
%
\def\eqnn#1{\xdef #1{(\secsym\the\meqno)}\writedef{#1\leftbracket#1}%
\global\advance\meqno by1\wrlabeL#1}
\def\eqna#1{\xdef #1##1{\hbox{$(\secsym\the\meqno##1)$}}
\writedef{#1\numbersign1\leftbracket#1{\numbersign1}}%
\global\advance\meqno by1\wrlabeL{#1$\{\}$}}
\def\eqn#1#2{\xdef #1{(\secsym\the\meqno)}\writedef{#1\leftbracket#1}%
\global\advance\meqno by1$$#2\eqno#1\eqlabeL#1$$}
%
\newskip\footskip\footskip14pt plus 1pt minus 1pt 
\def\footnotefont{\ninepoint}\def\f@t#1{\footnotefont #1\@foot}
\def\f@@t{\baselineskip\footskip\bgroup\footnotefont\aftergroup\@foot\let\next}
\setbox\strutbox=\hbox{\vrule height9.5pt depth4.5pt width0pt}
\global\newcount\ftno \global\ftno=0
\def\foot{\global\advance\ftno by1\footnote{$^{\the\ftno}$}}
%
\newwrite\ftfile
\def\footend{\def\foot{\global\advance\ftno by1\chardef\wfile=\ftfile
$^{\the\ftno}$\ifnum\ftno=1\immediate\openout\ftfile=foots.tmp\fi%
\immediate\write\ftfile{\noexpand\smallskip%
\noexpand\item{f\the\ftno:\ }\pctsign}\findarg}%
\def\footatend{\vfill\eject\immediate\closeout\ftfile{\parindent=20pt
\centerline{\bf Footnotes}\nobreak\bigskip\input foots.tmp }}}
\def\footatend{}
%
%
\global\newcount\refno \global\refno=1
\newwrite\rfile
%
\def\ref{\nref}
\def\nref#1{\xdef#1{[\the\refno]}\writedef{#1\leftbracket#1}%
\ifnum\refno=1\immediate\openout\rfile=refs.tmp\fi
\global\advance\refno by1\chardef\wfile=\rfile\immediate
\write\rfile{\noexpand\item{#1\ }\reflabeL{#1\hskip.31in}\pctsign}\findarg}
\def\findarg#1#{\begingroup\obeylines\newlinechar=`\^^M\pass@rg}
{\obeylines\gdef\pass@rg#1{\writ@line\relax #1^^M\hbox{}^^M}%
\gdef\writ@line#1^^M{\expandafter\toks0\expandafter{\striprel@x #1}%
\edef\next{\the\toks0}\ifx\next\em@rk\let\next=\endgroup\else\ifx\next\empty%
\else\immediate\write\wfile{\the\toks0}\fi\let\next=\writ@line\fi\next\relax}}
\def\striprel@x#1{} \def\em@rk{\hbox{}}
\def\lref{\begingroup\obeylines\lr@f}
\def\lr@f#1#2{\gdef#1{\ref#1{#2}}\endgroup\unskip}

\def\addref#1{\immediate\write\rfile{\noexpand\item{}#1}} 
\def\startrefs#1{\immediate\openout\rfile=refs.tmp\refno=#1}
\def\refs#1{\count255=1[\r@fs #1{\hbox{}}]}
\def\r@fs#1{\ifx\und@fined#1\message{reflabel \string#1 is undefined.}%
\nref#1{need to supply reference \string#1.}\fi%
\vphantom{\hphantom{#1}}\edef\next{#1}\ifx\next\em@rk\def\next{}%
\else\ifx\next#1\ifodd\count255\relax\xref#1\count255=0\fi%
\else#1\count255=1\fi\let\next=\r@fs\fi\next}
%

%
\newwrite\ffile\global\newcount\figno \global\figno=1
\def\fig{fig.~\the\figno\nfig}
\def\nfig#1{\xdef#1{fig.~\the\figno}%
\writedef{#1\leftbracket fig.\noexpand~\the\figno}%
\ifnum\figno=1\immediate\openout\ffile=figs.tmp\fi\chardef\wfile=\ffile%
\immediate\write\ffile{\noexpand\medskip\noexpand\item{Fig.\ \the\figno. }
\reflabeL{#1\hskip.55in}\pctsign}\global\advance\figno by1\findarg}
\def\xfig{\expandafter\xf@g}\def\xf@g fig.\penalty\@M\ {}
\def\figs#1{figs.~\f@gs #1{\hbox{}}}
\def\f@gs#1{\edef\next{#1}\ifx\next\em@rk\def\next{}\else
\ifx\next#1\xfig #1\else#1\fi\let\next=\f@gs\fi\next}
\newwrite\lfile
{\escapechar-1\xdef\pctsign{\string\%}\xdef\leftbracket{\string\{}
\xdef\rightbracket{\string\}}\xdef\numbersign{\string\#}}

\def\writestop{\def\writestoppt{\immediate\write\lfile{\string\pageno%
\the\pageno\string\startrefs\leftbracket\the\refno\rightbracket%
\string\def\string\secsym\leftbracket\secsym\rightbracket%
\string\secno\the\secno\string\meqno\the\meqno}\immediate\closeout\lfile}}
\def\writestoppt{}\def\writedef#1{}
\def\seclab#1{\xdef #1{\the\secno}\writedef{#1\leftbracket#1}\wrlabeL{#1=#1}}
\def\subseclab#1{\xdef #1{\secsym\the\subsecno}%
\writedef{#1\leftbracket#1}\wrlabeL{#1=#1}}
\newwrite\tfile \def\writetoca#1{}
\def\leaderfill{\leaders\hbox to 1em{\hss.\hss}\hfill}
\def\writetoc{\immediate\openout\tfile=toc.tmp
   \def\writetoca##1{{\edef\next{\write\tfile{\noindent ##1
   \string\leaderfill {\noexpand\number\pageno} \par}}\next}}}
%
%
%

\catcode`\@=12 
%
\edef\tfontsize{\ifx\answ\bigans scaled\magstep3\else scaled\magstep4\fi}
\font\titlerm=cmr10 \tfontsize \font\titlerms=cmr7 \tfontsize
\font\titlermss=cmr5 \tfontsize \font\titlei=cmmi10 \tfontsize
\font\titleis=cmmi7 \tfontsize \font\titleiss=cmmi5 \tfontsize
\font\titlesy=cmsy10 \tfontsize \font\titlesys=cmsy7 \tfontsize
\font\titlesyss=cmsy5 \tfontsize \font\titleit=cmti10 \tfontsize
\skewchar\titlei='177 \skewchar\titleis='177 \skewchar\titleiss='177
\skewchar\titlesy='60 \skewchar\titlesys='60 \skewchar\titlesyss='60
\def\titlefont{\def\rm{\fam0\titlerm}
\textfont0=\titlerm \scriptfont0=\titlerms \scriptscriptfont0=\titlermss
\textfont1=\titlei \scriptfont1=\titleis \scriptscriptfont1=\titleiss
\textfont2=\titlesy \scriptfont2=\titlesys \scriptscriptfont2=\titlesyss
\textfont\itfam=\titleit \def\it{\fam\itfam\titleit}\rm}
 \ifx\answ\bigans\else scaled\magstep1\fi
\ifx\answ\bigans\def\abstractfont{\tenpoint}\else
\font\abssl=cmsl10 scaled \magstep1
\font\absrm=cmr10 scaled\magstep1 \font\absrms=cmr7 scaled\magstep1
\font\absrmss=cmr5 scaled\magstep1 \font\absi=cmmi10 scaled\magstep1
\font\absis=cmmi7 scaled\magstep1 \font\absiss=cmmi5 scaled\magstep1
\font\abssy=cmsy10 scaled\magstep1 \font\abssys=cmsy7 scaled\magstep1
\font\abssyss=cmsy5 scaled\magstep1 \font\absbf=cmbx10 scaled\magstep1
\skewchar\absi='177 \skewchar\absis='177 \skewchar\absiss='177
\skewchar\abssy='60 \skewchar\abssys='60 \skewchar\abssyss='60
\def\abstractfont{\def\rm{\fam0\absrm}
\textfont0=\absrm \scriptfont0=\absrms \scriptscriptfont0=\absrmss
\textfont1=\absi \scriptfont1=\absis \scriptscriptfont1=\absiss
\textfont2=\abssy \scriptfont2=\abssys \scriptscriptfont2=\abssyss
\textfont\itfam=\bigit \def\it{\fam\itfam\bigit}\def\footnotefont{\tenpoint}%
\textfont\slfam=\abssl \def\sl{\fam\slfam\abssl}%
\textfont\bffam=\absbf \def\bf{\fam\bffam\absbf}\rm}\fi
\def\tenpoint{\def\rm{\fam0\tenrm}
\textfont0=\tenrm \scriptfont0=\sevenrm \scriptscriptfont0=\fiverm
\textfont1=\teni  \scriptfont1=\seveni  \scriptscriptfont1=\fivei
\textfont2=\tensy \scriptfont2=\sevensy \scriptscriptfont2=\fivesy
\textfont\itfam=\tenit \def\it{\fam\itfam\tenit}\def\footnotefont{\ninepoint}%
\textfont\bffam=\tenbf \def\bf{\fam\bffam\tenbf}\def\sl{\fam\slfam\tensl}\rm}
\font\ninerm=cmr9 \font\sixrm=cmr6 \font\ninei=cmmi9 \font\sixi=cmmi6
\font\ninesy=cmsy9 \font\sixsy=cmsy6 \font\ninebf=cmbx9
\font\nineit=cmti9 \font\ninesl=cmsl9 \skewchar\ninei='177
\skewchar\sixi='177 \skewchar\ninesy='60 \skewchar\sixsy='60
\def\ninepoint{\def\rm{\fam0\ninerm}
\textfont0=\ninerm \scriptfont0=\sixrm \scriptscriptfont0=\fiverm
\textfont1=\ninei \scriptfont1=\sixi \scriptscriptfont1=\fivei
\textfont2=\ninesy \scriptfont2=\sixsy \scriptscriptfont2=\fivesy
\textfont\itfam=\ninei \def\it{\fam\itfam\nineit}\def\sl{\fam\slfam\ninesl}%
\textfont\bffam=\ninebf \def\bf{\fam\bffam\ninebf}\rm}
%
%

\hyphenation{anom-aly anom-alies coun-ter-term coun-ter-terms}
\def\inv{^{\raise.15ex\hbox{${\scriptscriptstyle -}$}\kern-.05em 1}}

\def\Dsl{\,\raise.15ex\hbox{/}\mkern-13.5mu D} 
\def\dsl{\raise.15ex\hbox{/}\kern-.57em\partial}

\font\bigit=cmti10 scaled \magstep1
\def\lspace{\ifx\answ\bigans{}\else\qquad\fi}
\def\lbspace{\ifx\answ\bigans{}\else\hskip-.2in\fi} 
\def\boxeqn#1{\vcenter{\vbox{\hrule\hbox{\vrule\kern3pt\vbox{\kern3pt
    \hbox{${\displaystyle #1}$}\kern3pt}\kern3pt\vrule}\hrule}}}
\def\mbox#1#2{\vcenter{\hrule \hbox{\vrule height#2in
        \kern#1in \vrule} \hrule}}  
%

\def\darr#1{\raise1.5ex\hbox{$\leftrightarrow$}\mkern-16.5mu #1}

\def\half{{\textstyle{1\over2}}} 
\def\roughly#1{\raise.3ex\hbox{$#1$\kern-.75em\lower1ex\hbox{$\sim$}}}

%
%


\def\frac#1#2{{#1\over#2}}

\def\half{\frac12}

\def\journal#1&#2(#3){\unskip, #1~\bf #2 \rm(19#3) }
\def\andjournal#1&#2(#3){\sl #1~\bf #2 \rm (19#3) }

\def\bra#1{\left\langle #1\right|}
\def\ket#1{\left| #1\right\rangle}
\def\det{{\rm det}}

\catcode`\@=11\def\slash#1{\mathord{\mathpalette\c@ncel{#1}}}
\overfullrule=0pt
\def\steepslash{\c@ncel}
\def\frac#1#2{{#1\over #2}}

\def\:{\!:\!}
\def\inbar{\,\vrule height1.5ex width.4pt depth0pt}
\def\IQ{\relax\,\hbox{$\inbar\kern-.3em{\rm Q}$}}
\def\IB{\relax{\rm I\kern-.18em B}}
\def\IC{\relax\hbox{$\inbar\kern-.3em{\rm C}$}}
\def\IP{\relax{\rm I\kern-.18em P}}
\def\IR{\relax{\rm I\kern-.18em R}}
\def\ZZ{\relax\ifmmode\mathchoice
{\hbox{Z\kern-.4em Z}}{\hbox{Z\kern-.4em Z}}
{\lower.9pt\hbox{Z\kern-.4em Z}}
{\lower1.2pt\hbox{Z\kern-.4em Z}}\else{Z\kern-.4em Z}\fi}

\catcode`\@=12

\def\npb#1(#2)#3{{ Nucl. Phys. }{B#1} (#2) #3}
\def\plb#1(#2)#3{{ Phys. Lett. }{#1B} (#2) #3}
\def\pla#1(#2)#3{{ Phys. Lett. }{#1A} (#2) #3}
\def\prl#1(#2)#3{{ Phys. Rev. Lett. }{#1} (#2) #3}
\def\mpla#1(#2)#3{{ Mod. Phys. Lett. }{A#1} (#2) #3}
\def\ijmpa#1(#2)#3{{ Int. J. Mod. Phys. }{A#1} (#2) #3}
\def\cmp#1(#2)#3{{ Comm. Math. Phys. }{#1} (#2) #3}
\def\cqg#1(#2)#3{{ Class. Quantum Grav. }{#1} (#2) #3}
\def\jmp#1(#2)#3{{ J. Math. Phys. }{#1} (#2) #3}
\def\anp#1(#2)#3{{ Ann. Phys. }{#1} (#2) #3}
\def\prd#1(#2)#3{{ Phys. Rev. } {D{#1}} (#2) #3}
\def\ptp#1(#2)#3{{ Progr. Theor. Phys. }{#1} (#2) #3}
\def\aom#1(#2)#3{{ Ann. Math. }{#1} (#2) #3}

\def\noi{\noindent}
\def\bs{\bigskip}
\def\bu{\bs\noindent $\bullet$}
\def\br{\buildrel}
\def\bra{\langle}
\def\ket{\rangle}

\def\G{{\bf G}}

\def\Q{{\bf Q}}
\def\R{{\bf R}}

\def\Z{{\bf Z}}

\def\cI{{\cal I}}

\def\cL{{\cal L}}

\def\cN{{\cal N}}

\def\cP{{\cal P}}

\def\cR{{\cal R}}

\def\cT{{\cal T}}

\def\cZ{{\cal Z}}
\input amssym

\def\cicy#1(#2|#3)#4{\left(\matrix{#2}\right|\!\!
                     \left|\matrix{#3}\right)^{{#4}}_{#1}}

\def\lra{\longrightarrow}

\def\hra{\hookrightarrow}

\def\ra{\rightarrow}

\def\bs{\bigskip}

\def\Box{{\,\lower0.9pt\vbox{\hrule
\hbox{\vrule height 0.2 cm \hskip 0.2 cm
\vrule height 0.2 cm}\hrule}\,}}

\global\newcount\thmno \global\thmno=0
\def\definition#1{\global\advance\thmno by1
\bigskip\noindent{\bf Definition \secsym\the\thmno. }{\it #1}
\par\nobreak\medskip\nobreak}
\def\question#1{\global\advance\thmno by1
\bigskip\noindent{\bf Question \secsym\the\thmno. }{\it #1}
\par\nobreak\medskip\nobreak}
\def\theorem#1{\global\advance\thmno by1
\bigskip\noindent{\bf Theorem \secsym\the\thmno. }{\it #1}
\par\nobreak\medskip\nobreak}
\def\proposition#1{\global\advance\thmno by1
\bigskip\noindent{\bf Proposition \secsym\the\thmno. }{\it #1}
\par\nobreak\medskip\nobreak}
\def\corollary#1{\global\advance\thmno by1
\bigskip\noindent{\bf Corollary \secsym\the\thmno. }{\it #1}
\par\nobreak\medskip\nobreak}
\def\lemma#1{\global\advance\thmno by1
\bigskip\noindent{\bf Lemma \secsym\the\thmno. }{\it #1}
\par\nobreak\medskip\nobreak}
\def\conjecture#1{\global\advance\thmno by1
\bigskip\noindent{\bf Conjecture \secsym\the\thmno. }{\it #1}
\par\nobreak\medskip\nobreak}
\def\exercise#1{\global\advance\thmno by1
\bigskip\noindent{\bf Exercise \secsym\the\thmno. }{\it #1}
\par\nobreak\medskip\nobreak}
\def\remark#1{\global\advance\thmno by1
\bigskip\noindent{\bf Remark \secsym\the\thmno. }{\it #1}
\par\nobreak\medskip\nobreak}
\def\problem#1{\global\advance\thmno by1
\bigskip\noindent{\bf Problem \secsym\the\thmno. }{\it #1}
\par\nobreak\medskip\nobreak}
\def\others#1#2{\global\advance\thmno by1
\bigskip\noindent{\bf #1 \secsym\the\thmno. }{\it #2}
\par\nobreak\medskip\nobreak}
\def\proof{\noindent Proof: }

\def\thmlab#1{\xdef #1{\secsym\the\thmno}\writedef{#1\leftbracket#1}\wrlabeL{#1=#1}}
%
%
\def\newsec#1{\global\advance\secno by1\message{(\the\secno. #1)}
\global\subsecno=0\thmno=0\eqnres@t\noindent{\bf\the\secno. #1}
\writetoca{{\secsym} {#1}}\par\nobreak\medskip\nobreak}
\def\eqnres@t{\xdef\secsym{\the\secno.}\global\meqno=1\bigbreak\bigskip}
\def\sequentialequations{\def\eqnres@t{\bigbreak}}\xdef\secsym{}
%

%
\newcount{\exnum}
\def\prob{\advance\exnum by 1
\bigskip\item{\the\exnum.}\ }

\baselineskip=12pt plus 2pt minus 2pt
\parskip 5pt
\def\ol{\overline}


\Title{}{Counting Unimodular Lattices in $\R^{r,s}$}

\centerline{
Shinobu Hosono, Bong H. Lian, Keiji Oguiso, and Shing-Tung Yau
}

\vskip.2in
\noindent Abstract.
Narain lattices are unimodular lattices {\it in} $\R^{r,s}$,
subject to certain natural equivalence relation
and rationality condition.
The problem of describing and
counting these rational equivalence classes of
Narain lattices in
$\R^{2,2}$ has led to an interesting connection to binary forms
and their Gauss products, as shown in [HLOYII]. As a sequel, in
this paper, we study arbitrary rational Narain lattices and
generalize some of our earlier results. In particular in the case
of $\R^{2,2}$, a new interpretation of the Gauss product of binary
forms brings new light to a number of related objects -- rank 4
rational Narain lattices, over-lattices, rank 2 primitive
sublattices of an abstract rank 4 even unimodular lattice $U^2$,
and isomorphisms of discriminant groups of rank 2 lattices.

\vskip.2in

\newsec{Introduction}

{\it The Main Problem.}
Let $E$ a real quadratic space, i.e. a real vector space equipped
with a non-degenerate bilinear form $\bra,\ket$ of signature $(r,s)$.
Assume that $8$ divides $r-s$. Fix a non-degenerate linear subspace $V\subset E$.
An even integral unimodular lattice $\Gamma\subset E$ of signature $(r,s)$
is called a Narain lattice in $E$.
If $rank~\Gamma\cap V=dim~V$, we say that $\Gamma$ is $V$-rational.
The problem is to describe and count
$V$-rational Narain lattices up to the action of
the orthogonal subgroup $O(V)\times O(V^\perp)\subset O(E)$.
For $r,s>0$, Narain lattices are parameterized by the set
$O(E)/O(\Gamma_0)$, where $\Gamma_0$ is a fixed Narain lattice.
Geometrically in this case, the problem amounts to
describing and counting the ``rational points'' in the homogeneous space
$O(V)\times O(V^\perp)\backslash O(E)/O(\Gamma_0)$, where $\Gamma_0$
is a fixed $V$-rational Narain lattice.

For simplicity, we will concentrate on the case $E=\R^{n,n}$
throughout most of the paper, and return to the general case only
at the end.
The study of general $V$-rational Narain lattices in $\R^{n,n}$
leads us to consider the following classes of objects,
each of which comes equipped with its own suitable equivalence relation
(to be described later):

\item{i.} Rational Narain lattices in $\R^{n,n}$;
\item{ii.} Rank $2n$ even unimodular over-lattices and triples;
\item{iii.} Rank $n$ primitive sublattices of $U^n$;
\item{iv.} Isomorphisms of discriminant groups;
\item{v.} Pairs of binary quadratic forms and their Gauss products, when $n=2$.

\nopagenumbers
\headline{\ifodd\pageno\rightheadline\else\leftheadline\fi}
\def\rightheadline{\tenrm\hfil Counting Rational Narain Lattices\hfil\folio}
\def\leftheadline{\tenrm\folio\hfil S. Hosono, B. Lian, K. Oguiso, S.T. Yau \hfil}

\noi Here $U$ is an abstract even integral unimodular
hyperbolic rank 2 lattice.
Narain lattices in $E=\R^{n,n}$
arouse in physics as a way to build modular invariant
string theory models [Na][NSW][Po]
and conformal field theories [Mo][GV][W]. Here one
sets $V=\R^{n,0}\subset E$, and let $E\ra V$, $x\mapsto x_L$,
and $E\ra V^\perp$, $x\mapsto x_R$, denote the orthogonal projections.
In string theory [Na],
to each even unimodular lattice $\Gamma$ {\it in} $E$,
one attaches the real analytic function
$$
Z_\Gamma(\tau)={1\over|\eta(\tau)|^{2n}}\sum_{x\in\Gamma}
e^{\pi i\tau\bra x_L,x_L\ket} e^{\pi i\bar\tau\bra x_R,x_R\ket}
$$
known as a partition function, defined on the complex upper half plane.
By Poisson summation (cf. [Po][S]), one shows that
$Z_\Gamma(\tau)$ is modular invariant. An important problem in
string theory is to understand the behavior
of the modular functions $Z_\Gamma(\tau)$
as one deforms the ``moduli'' $\Gamma$ (see [HLOYII] and
references therein).

The theory of over-lattices, also known as gluing theory,
came up in coding theory (see [CS] p99 and references therein),
and in the study of primitive embeddings of lattices [Ni],
in which discriminant groups also play an important role.
The latter was a key ingredient
in our recent solution to the counting problem
of Fourier Mukai partners for K3 surfaces
in [HLOYI]. More recently, the same objects (all but iii.)
have also arouse in the study
of rational conformal field theory on real tori [HLOYII].
There, the Gauss product of positive definite coprime binary forms
were essential in giving a geometric description of those
rational conformal field theories.

\def\gammaTheorem{3.12}
\def\thetaTheorem{4.8}
\def\LambdaGauss{5.13}
\def\TwoOneTheorem{6.5}
\def\LambdaTheorem{6.6}

The approach taken in this paper is entirely algebraic.
Let's begin with an idea that goes back essentially to Gauss.
Given two binary quadratic forms $f(x,y),g(x,y)$ of the same discriminant,
Gauss defined his composition law by constructing a third form $h(x,y)$
such that
$$
f(x,y)g(x',y')=h(X,Y)
$$
where $X,Y$ are certain integral linear forms of $x,y,x',y'$.
This classical construction has, of course, long been
superseded by the modern theory of ideal class groups
for general Dedekind fields. The key new idea in this paper is this:
that Gauss' original construction of his composition law, in fact, takes on
a particularly interesting and potent
form, when interpreted
{\it inside} the rank 4 lattice of $2\times 2$ integral
matrices with the quadratic form $2det~X$. This special rank 4 lattice
is abstractly isomorphic to $U^2$. But because it is
also a {\it ring} with lots of symmetry, all the algebraic
structures that come with it can be brought to bear on the
study of sublattices. Moreover there is a very interesting {\it duality}:
a binary form can be viewed
as a point (an integral symmetric matrix), and a rank 2 primitive
sublattice (equipped with a basis) can be view as a quadratic form.
It is this duality, together with the algebraic structures of the matrix ring,
that makes our algebraic approach work.
The details are spelled out in Section 5.


Here is an outline of the paper.
Theorem \gammaTheorem~
gives a 1-1 correspondence between equivalence classes
of objects in i. and those in ii. Theorem \thetaTheorem~
does that same for objects in ii. and iii. Theorem \LambdaGauss~
shows precisely how the Gauss products
of coprime binary forms enter the description of iii.
through a new map $\Lambda$ when $n=2$. This is where
the quadratic form $2det~X$ on matrices comes in.
We also describe briefly the connection between ii. and iv.
Finally, Theorems \TwoOneTheorem~ and \LambdaTheorem~ use $\Lambda$ to
give a description of iii., hence culminating in
a description of $V$-rational Narain lattices in $\R^{2,2}$,
in terms of binary forms.
In the last section, we apply this to the positive
definite case and recover a result in [HLOYII].
We also discuss the indefinite case as a new application.
Finally, we comment on some further interesting generalizations
in the last section.

{\it Acknowledgement.}
We thank B. Gross for some helpful suggestions.
The third named authors thank the
the Education Ministry of Japan
for financial support. The second and the fourth named authors
are supported by NSF grants.

\subsec{General conventions}

\bu Let $X,Y$ be sets equipped with a transformation group
$G$. Let $f:X\ra Y$ be a map of sets. We say that
$f$ {\it descends through} $G$ if $f$ sends each $G$-orbit into a $G$-orbit,
i.e. $\forall g\in G$, $\forall x\in X$, $f(g\cdot x)=g'\cdot f(x)$ for some $g'\in G$.
In this case, the composition $X\ra Y\ra Y/G$
descends to an induced map $X/G\ra Y/G$. We say that $f$ is
{\it $G$-equivariant} if $f$ commutes with the $G$-action,
i.e. $\forall g\in G$, $\forall x\in X$, $f(g\cdot x)=g\cdot f(x)$. In
this case, $f$ descends through {\it any} subgroup $K\subset G$.
We say that $f$ {\it factors through} $G$ if $f$ sends each $G$-orbit to a point,
i.e. $\forall g\in G$, $\forall x\in X$, $f(g\cdot x)=f(x)$.
In this case, $f$ descends to an induced map $X/G\ra Y$.

\bu If the transformation group $G$ on $X$ is of the form $H\times K$,
we can let $H$ acts on the left and $K$ on the right, and denote
a $G$-orbit by $H\cdot x\cdot K$. The orbit space in this case is denoted by
the double quotient $H\backslash X/K$. Likewise, we can also speak of
the orbit spaces $H\backslash X$ and $X/K$. If $X$ is a group
and $H,K$ are subgroups, then these are the usual left and right coset
spaces.

\bu If $f:X\ra Y$ is an isomorphism of lattices, we denote by
$f^*:X^*\ra Y^*$ (the direction is not reversed!) the inverse
of the dual isomorphism between their dual lattices.
We denote by $f^*:A_X\ra A_Y$ the
the induced isometry of discriminant groups.

\bu We will only consider {\it even} integral lattices $X$,
i.e. a free abelian group of finite rank equipped
with an even {\it non-degenerate} integral quadratic form $\bra,\ket$.
The discriminant of a lattice $X$ is defined to be $(-1)^{n-1}$ times
the determinant of the matrix $\bra v_i,v_j\ket$ for
a given $\Z$-base $v_i$ of $X$. The sign is chosen to
make it consistent with the rank 2 case. It is clear that the discriminant
is invariant under a change of base, hence it is
a well-defined invariant of the lattice.
If $X$ is an even lattice equipped with $\bra,\ket$,
and $r\in\Q$ is a nonzero number such that
the scalar multiple $r\bra,\ket$
remains even integral, then
we denote by $X(r)$ the same
abelian group but equipped with the form $r\bra,\ket$.
In particular, $X(-1)$ is the lattice with
the sign of the form reversed.
We say that $X$ is {\it primitive} if
$m=1$ is the only positive integer
such that $X({1\over m})$ is even integral. Thus
the integral binary quadratic form
$[a,b,c]\equiv\left(\right.\Z^2,\left[\matrix{2a&b\cr b&2c}\right]\left.\right)$
is primitive iff $gcd(a,b,c)=1$.

\bu The term primitive is used in one other (somewhat confusing) way.
We say that $X\subset Y$ is a {\it primitive sublattice} of $Y$
if the inclusion is isometric and $Y/X$ is torsion free.
Note that this does not mean that $Y$ is primitive,
as an abstract lattice.

\bu The symbol $\sigma$ is used throughout the paper, but
to mean different things in a few different but related contexts.
For example, $\sigma:\R^{n,n}\ra\R^{n,n}$
will be an involutive anti-isometry. It induces an
involution on the set of rational Narain lattices
which we denote by $\sigma:\cR\ra\cR$. On the set $\cT$ of triples
$(X,Y,\varphi)$, of lattices $X,Y$ and isometry $\varphi:A_Y\ra A_X$,
we have an involution $(X,Y,\varphi)\mapsto(Y,X,\varphi^{-1})$.
We denote this by $\sigma:\cT\ra\cT$. The permutation
$(P,Q)\mapsto(Q,P)$ acting on the set $\cP$ of pairs of quadratic forms.
We denote this by $\sigma:\cP\ra\cP$. All of these will eventually
be related.

\newsec{Narain Lattices}

\bu Notations:

$\R^{r,s}$: $\R^{r+s}$ equipped with the quadratic form
$\bra x,x\ket=x_1^2+\cdots+x_r^2-x_{r+1}^2-\cdots-x_{r+s}^2$.

$U$: the hyperbolic even unimodular lattice $\Z e\oplus\Z f$ with $\bra e,e\ket=\bra f,f\ket=0$,
$\bra e,f\ket=1$.

$U^n$: the direct sum of $n$ copies of $U$;
we name the $i$th copy $\Z e_i\oplus\Z f_i$.


$\varepsilon_1,..,\varepsilon_{2n}$: the standard basis of $\R^{2n}$;
$\varepsilon_i=(0,..,0,1,0,..,0)$, 1 being at the $i$th slot.

$E_i^e:={1\over\sqrt2}(\varepsilon_i+\varepsilon_{n+i}),~~~
F_i^e:={1\over\sqrt2}(\varepsilon_i-\varepsilon_{n+i})$: called the standard Narain basis of $\R^{n,n}$.

$O(\R^{n,n})$: the group of linear isometries of $\R^{n,n}$.


$O(\Gamma^e)$: the subgroup of $O(\R^{n,n})$ preserving the special even unimodular lattice
$$
\Gamma^e:=\sum_i(\Z E_i^e+\Z F_i^e)\subset\R^{n,n}.
$$

\definition{A Narain embedding is an isometric embedding
$$
\nu:U^n\hra\R^{n,n}.
$$
The set of Narain embeddings is denoted by $\cI_n$.
An even unimodular
lattice $\Gamma\subset\R^{n,n}$ of signature $(n,n)$
is called a Narain lattice. The set of Narain lattices is
denoted by $\cN_n$. The subscript $n$ will be surpressed
when there is no confusion.}

\bu Since a Narain embedding $\nu$ is an isometry, the images
$\nu(e_i),\nu(f_i)\in\R^{n,n}$ form a basis having the same
inner products as the basis vectors $e_i,f_i$ of $U^n$. We shall call
such an {\it ordered $\R$-basis} of $\R^{n,n}$
{\it a Narain basis}.
Conversely, any Narain basis $E_i,F_i$ of $\R^{n,n}$
defines a Narain embedding $\nu$ by declaring
$$
\nu(e_i)=E_i,~~~\nu(f_i)=F_i.
$$
Thus {\it Narain embeddings correspond 1-1 with Narain bases.}
Denote by $\nu_e$ the Narain embedding corresponding to
the special Narain basis $E_i^e, F_i^e$.

\lemma{There is a unique $O(\R^{n,n})$-equivariant bijection
$O(\R^{n,n})\ra\cI$ such that $e\mapsto\nu_e$.}
\proof
An element $g\in O(\R^{n,n})$ acts on
Narain embeddings $\nu$ by left translations:
$$
\nu\mapsto g\circ\nu.
$$
The action is transitive because any two Narain bases
of $\R^{n,n}$ are related by a {\it unique} $g\in O(\R^{n,n})$.
Conversely, given a Narain basis, its image under a $g\in O(\R^{n,n})$ forms
another Narain basis.
Thus the map
$$
O(\R^{n,n})\ra\cI,~~~g\mapsto\nu_g:=g\circ\nu_e.
$$
is a bijection. It clearly has the asserted uniqueness and equivariance property. $\Box$

\corollary{The composition $O(\R^{n,n})\ra\cI\ra\cN$, $g\mapsto\nu_g\mapsto\nu_g(U^n)$,
descends to a bijection
$$
O(\R^{n,n})/O(\Gamma^e)\br\sim\over\ra\cN,~~~g\cdot O(\Gamma^e)\mapsto\nu_g(U^n).
$$}
\proof
It is clear that the map $O(\R^{n,n})\ra\cN$ factors through $O(\Gamma^e)$, which
preserves the lattice $\nu_e(U^n)=\Gamma^e$. By Milnor's theorem
every Narain lattice $\Gamma\subset\R^{n,n}$ is in the image of $\cI\ra\cN$.
Thus $O(\R^{n,n})/O(\Gamma^e)\ra\cN$ is surjective. Equivalently
$O(\R^{n,n})$ acts transitively on $\cN$. With the based point $\nu_e(U^n)\in\cN$,
the isotropy group is $O(\Gamma^e)$. Thus we have a bijection. $\Box$

\subsec{$V$-equivalence}

\bu Once and for all, we fix a {\it non-degenerate}
linear subspace $V\subset\R^{n,n}$, i.e. $\bra,\ket|V$ non-degenerate. The notation
$V^\perp$ will always mean the orthogonal complement of $V$ in $\R^{n,n}$.
At the end, for $n=2$,
we shall specialize $V$ to the case $V=\R^{2,0}=\{(*,*,0,0)\}\subset\R^{2,2}$
or $V=\R^{1,1}=\{(*,0,*,0)\}\subset\R^{2,2}$.

\bu Since $\R^{n,n}=V\oplus V^\perp$ canonically, we have the canonical inclusions
of isometry groups
$$
O(V),~O(V^\perp)\subset O(V)\times O(V^\perp)\subset O(\R^{n,n}).
$$
Note that the middle group is the subgroup of $O(\R^{n,n})$
preserving the decomposition $V\oplus V^\perp$, i.e. an
element $f\in O(\R^{n,n})$ is in $O(V)\times O(V^\perp)$ iff
$fV=V$ and $fV^\perp=V^\perp$.

\definition{Two Narain lattices $\Gamma,\Gamma'$
are called $V$-equivalent if $\Gamma'=g\Gamma$ for some $g\in O(V)\times O(V^\perp)$.}

\corollary{The map $O(\R^{n,n})/O(\Gamma^e)\br\sim\over\ra\cN$ above
 descends to $V$-equivalence classes, i.e.
$$\eqalign{
&O(V)\times O(V^\perp)\backslash O(\R^{n,n})/O(\Gamma^e)\br\sim\over\ra O(V)\times O(V^\perp)\backslash\cN,\cr
&O(V)\times O(V^\perp)\cdot g\cdot O(\Gamma^e)\mapsto O(V)\times O(V^\perp)\cdot\nu_g(U^n).
}$$}

\bu Whenever convenient,  we shall identify Narain bases
with Narain embeddings,  and with $O(\R^{n,n})$ via
$$
(\nu_g(e_1),..,\nu_g(e_n),\nu_g(f_1),..,\nu_g(f_n))\equiv\nu_g\equiv g,
$$
and $\cN$ with $O(\R^{n,n})/O(\Gamma^e)$ via
$$
\nu_g(U^n)\equiv g\cdot O(\Gamma^e).
$$

\subsec{An involutive anti-isometry $\sigma$}

Throughout this paper, we assume that $V$ has
$dim~V=n$ and has signature $(p,n-p)$. Then $V^\perp$ has signature $(n-p,p)$,
hence $V\cong V^\perp(-1)$ as quadratic spaces.
{\it We fix an anti-isometry}
$\sigma:V\ra V^\perp$, and define
the linear map
$$
\R^{n,n}\ra\R^{n,n},~~\R^{n,n}=V\oplus V^\perp\ni(x,y)\mapsto(\sigma^{-1}y,\sigma x).
$$
This is clearly an involutive anti-isometry of $\R^{n,n}$
exchanging the two subspaces $V,V^\perp\subset\R^{n,n}$.
We also denote this involution by $\sigma$.
Note that if $\Gamma$ is a Narain lattice, then
its image $\sigma\Gamma\subset\R^{n,n}$
is again an even unimodular lattice as a subset
of the quadratic space $\R^{n,n}$, hence $\sigma\Gamma$
is a Narain lattice.
So we have an involution $\cN\ra\cN$,
$\Gamma\mapsto\sigma\Gamma$, which we also denote by $\sigma$.

\lemma{The correspondence $g\mapsto \sigma g\sigma$,
is a well-defined involution on $O(\R^{n,n})$,
hence on the set of Narain embeddings. This involution
also stabilizes the subgroups $O(V)\times O(V^\perp)$, and
the identity component $O_0(\R^{n,n})$. If $\sigma\Gamma^e=\Gamma^e$, then
$\sigma$ stabilizes $O(\Gamma^e)$ as well.}
\thmlab\NormalizationLemma
\proof
For $x\in\R^{n,n}$ and $g\in O(\R^{n,n})$, we have
$$
\bra \sigma g\sigma x,\sigma g\sigma x\ket=
-\bra g\sigma x,g\sigma x\ket=-\bra\sigma x,\sigma x\ket=\bra x,x\ket.
$$
This shows that $\sigma g\sigma\in O(\R^{n,n})$.
That conjugation by $\sigma$ stabilizes $O(V)\times O(V^\perp)$
is obvious from the definition of $\sigma$.
Since $\sigma\cdot O_0(\R^{n,n})\cdot\sigma$ is a connected
subgroup of $O(\R^{n,n})$ having the same dimension,
this subgroup must coincide with $O_0(\R^{n,n})$.
The last statement is clear. $\Box$

\newsec{$V$-Rational Narain Lattices and Triples}

Recall that we have fixed a non-degenerate subspace $V\subset\R^{n,n}$.

\definition{Given a Narain lattice $\Gamma$, the primitive sublattices
$$
\Gamma_L:=\Gamma\cap V,~~~\Gamma_R:=\Gamma\cap V^\perp
$$
are respectively called the left and right $V$-lattices of $\Gamma$.
A Narain lattice $\Gamma$ is called $V$-rational if
$\Gamma_L$ has maximal rank, i.e. $rk~\Gamma_L=dim~V$.
The set of $V$-rational Narain lattices is denoted by
$$
\cR\subset\cN\equiv O(\R^{n,n})/O(\Gamma^e).
$$
The discriminant of a $V$-rational Narain lattice $\Gamma$
is defined to be the discriminant of the lattice $\Gamma_L$.
In this case, we will say that $\Gamma$ is primitive
if $\Gamma_L$ is primitive as an abstract lattice.
Note that we have suppressed $V$ from the notation, even though $\cR$
depends on $V$. We will also drop the mention of $V$
when there is no confusion.
}

If $V$ has signature $(p,q)$,
and if $\Gamma$ is $V$-rational, then $\Gamma_L$
also has signature $(p,q)$, by definition.
Similarly $V^\perp$ and $\Gamma_R=\Gamma_L^\perp$ in $\Gamma$
both have signature $(n-p,n-q)$. It follows that
a $V$-rational Narain lattice $\Gamma$ is an {\it even unimodular overlattice}
of the lattice $\Gamma_L\oplus\Gamma_R\subset\Gamma$.

\lemma{If $V,V'\subset\R^{n,n}$ are non-degenerate subspaces with
the same signature, there is an isometry $f\in O(\R^{n,n})$
such that $fV=V'$, and that $f\Gamma$ is a $V'$-rational
for every $V$-rational Narain lattice $\Gamma$.}
\thmlab\Vindependence
\proof
Let $(p,q)$ be the signature of $V$,
and let $e_1,...,e_{2n}$ be an orthogonal basis for $\R^{n,n}$ such that
$e_i\in V$ for $1\leq i\leq p+q$, that $e_i\in V^\perp$ for $p+q+1\leq i\leq 2n$,
and that
$$
\bra e_i,e_i\ket=\left\{\matrix{
+1 & 1\leq i\leq p\cr
-1 & p+1\leq i\leq p+q\cr
+1 & p+q+1\leq i\leq n+q\cr
-1 & n+q+1\leq i\leq 2n.}\right.
$$
Likewise for $e_1',...,e_{2n}'$ and $V'$.
Then $f:e_i\mapsto e_i'$ defines an element $f\in O(\R^{n,n})$
with $fV=V'$.

If $\Gamma$ is a $V$-rational Narain lattice, then
$\Gamma\cap V$ has maximal rank, hence
$f\Gamma\cap fV=f\Gamma\cap V'$ also has maximal rank,
which means that $f\Gamma$ is $V'$-rational. $\Box$

\lemma{For each $\Gamma\in\cR$, there exists a unique isometry
$\varphi:A_{\Gamma_R(-1)}\ra A_{\Gamma_L}$ such that
$$
(*)~~~~\Gamma/(\Gamma_L\oplus\Gamma_R)=\{\varphi(a)\oplus a|a\in \Gamma_R(-1)^*/\Gamma_R(-1)\}.
$$
In particular, we have
$$
\det~\Gamma_L=|A_{\Gamma_L}|=|A_{\Gamma_R(-1)}|=det~\Gamma_R(-1).
$$
Moreover, when $dim~V=n$, then $\Gamma_L, \Gamma_R(-1)$ have the same signature,
and they are in the same genus. In this case,
$\Gamma_L$ is primitive iff $\Gamma_R$ is primitive.}
\thmlab\Overlattice
\proof
The first assertion follows from [Ni]:
even unimodular overlattices correspond 1-1 with
isometries $\varphi:A_{\Gamma_R(-1)}\br\sim\over\ra A_{\Gamma_L}$
such that (*) holds. The second assertion follows the characterization [Ni]
of genus that two lattices $X,Y$ are in the same genus iff
they have the same signature and have $A_X\cong A_Y$. The last assertion is
a result of the next lemma. $\Box$

\lemma{Two lattices $X,Y$ in the same genus
are either both primitive or both not primitive.}
\thmlab\BothPrimitive
\proof
Suppose that $X,Y$ are in the same genus, and $X$ is primitive, but $Y$ is not.
By definition, we have an isomorphism $\varphi_p:X\otimes\Z_p\ra Y\otimes\Z_p$
of $\Z_p$ integral lattices, for each prime number $p$.
Let $q$ be a prime number such that $Y({1\over q})$ remains even integral,
and let $u\in X$ be a vector such that ${1\over q}\bra u,u\ket_X\notin2\Z$.
First consider the case $q\neq2$.
Then ${1\over q}\bra u,u\ket_X\otimes\Z_q\notin2\Z_q$.
(Otherwise, the integer number $\bra u,u\ket_X$ would have the ``zero'' to cancel
the ``pole'' ${1\over q}$ in $\Q_q$.) On the other hand ${1\over q}\bra,\ket_Y$
is an even integral form. In particular, we have ${1\over q}\bra v,v\ket_Y\otimes\Z_q\in2\Z_q$
for any $v\in Y\otimes\Z_q$. Since $\varphi_q(u)\in Y\otimes\Z_q$, we have
$$
{1\over q}\bra u,u\ket_X\otimes\Z_q={1\over q}\bra\varphi_q(u),\varphi_q(u)\ket_Y\otimes\Z_q\in 2\Z_q
$$
which is a contradiction.

For $q=2$, the same argument works, with the modification that
 now $\bra u,u\ket_X$
would have had to have ``zero'' of at least order 2 to cancel the
``pole'' ${1\over2}$. $\Box$

Recall that if $dim~V=n$ and has signature $(p,n-p)$,
and we fix an anti-isometry
$\sigma:V\ra V^\perp$, then we get an involutive anti-isometry
$\sigma:\R^{n,n}\ra\R^{n,n}$,
$(x,y)\mapsto(\sigma^{-1}y,\sigma x)$,
where $(x,y)\in V\oplus V^\perp=\R^{n,n}$.

\lemma{The action of $O(V)\times O(V^\perp)$ on the set $\cN$
preserves $V$-rationality, discriminant of the left and right $V$-lattices,
and primitivity of the left and right $V$-lattices. Likewise
the involution $\sigma:\cN\ra\cN$ preserves all these properties, except for discriminant,
which changes by the sign $(-1)^n$.
Moreover the involution descends to $V$-equivalence classes.}
\thmlab\sigmaLemma
\proof
Let $b\in O(V)\times O(V^\perp)$.
Since $b$ preserves the orthogonal decomposition $V\oplus V^\perp$,
it follows that
$$\eqalign{
&(b\Gamma)_L:=(b\Gamma)\cap V=b(\Gamma\cap V)=b\Gamma_L\cr
&(b\Gamma)_R:=(b\Gamma)\cap V^\perp=b(\Gamma\cap V^\perp)=b\Gamma_R.
}$$
Since
$$
rk~b\Gamma_L=rk~\Gamma_L,~~~det~b\Gamma_L=det~\Gamma_L,
$$
it follows that $b$ preserves rationality, discriminant of the left
lattices; likewise for the right lattices.
By definition, the quadratic form on $(b\Gamma)_L$ is
nothing but $\bra b^{-1}-,b^{-1}-\ket$ where $\bra,\ket$ is the quadratic form on $\Gamma_L$.
This shows that $\Gamma_L$ is primitive iff $(b\Gamma)_L$ is primitive.

Suppose $dim~V=n$, and $\Gamma$ is rational. Then
$$
\Gamma_L=\Gamma\cap V
$$
has maximal rank. Applying $\sigma$ to this, we get
$$
\sigma\Gamma_L=\sigma\Gamma\cap V^\perp =(\sigma\Gamma)_R.
$$
This shows that the right (and left) lattice of $\sigma\Gamma$
have rank $n$. Since $\sigma$ merely reverses the overall sign
of the quadratic form on $\R^{n,n}$, i.e. $\bra\sigma x,\sigma y\ket=\bra x,y\ket$,
primitivity of left and right lattices are preserved by $\sigma$,
and the discriminant changes by the sign $(-1)^n$.

By Lemma \NormalizationLemma, the element
$\sigma\in GL(2n,\R)$ normalizes the subgroup $O(V)\times O(V^\perp)$.
Thus the action of $\sigma$ on $\cN$ descends to the $O(V)\times O(V^\perp)$
orbits of Narain lattices, i.e. to $V$-
equivalence classes. $\Box$

\definition{
The set of $V$-equivalence classes of rational Narain
lattices is denoted by $\ol\cR\subset\ol\cN$.}

\definition{A triple $(X,Y,\varphi)$ consists of
even lattices $X,Y$ having the same signature and rank $n$,
and equipped with an isometry $\varphi:A_Y\ra A_X$.
The set of triples is denoted by $\cT$.
A triple $(X,Y,\varphi)$ is said to be primitive
if $X,Y$ are primitive.
Two triples $(X,Y,\varphi),(X',Y',\varphi')$ are said to
be properly equivalent if there exist isomorphisms $f:X\ra X', g:Y\ra Y'$
such that
$$
\varphi'=f^*\circ\varphi\circ {g^*}^{-1}.
$$
The set of proper equivalence classes of triples is denoted by $\ol\cT$.
The discriminant of a triple $(X,Y,\varphi)$ is defined to be the discriminant
of the lattices $X,Y$.}

\definition{Define the involution
$$
\sigma:\cT\ra\cT,~~~(X,Y,\varphi)\mapsto(Y,X,\varphi^{-1})
$$
Likewise on $\ol\cT$.
The triple $(X',Y',\varphi')$ is said to be
(improperly) equivalent to $(X,Y,\varphi)$ if it
is properly equivalent to either $(X,Y,\varphi)$ or $\sigma(X,Y,\varphi)$.}

\subsec{The map $\gamma$ and its invariance properties}


\definition{Define the map $\gamma:\cR\ra\cT$,
$\Gamma\mapsto(\Gamma_L,\Gamma_R,\varphi)$
where $\varphi:A_{\Gamma_R(-1)}\ra A_{\Gamma_L}$ is
the isometry determined by the overlattice $\Gamma\supset\Gamma_L\oplus\Gamma_R$.
}\thmlab\gammaDefinition

By Lemma \BothPrimitive, it follows that $\gamma(\Gamma)$ is primitive
iff $\Gamma_L$ (or $\Gamma_R$) is primitive.

We now consider how the map $\gamma$ interacts with
the involution $\sigma:\cR\ra\cR$ and
the group action of $O(V)\times O(V^\perp)$ on $\cR$.

\lemma{
\item{i.} If $\Gamma\in\cR$ and $b\in O(V)\times O(V^\perp)$, then
$\gamma(b\Gamma)$ is properly equivalent to $\gamma(\Gamma)$.
\item{ii.} If $\Gamma\in\cR$, then
$\gamma(\sigma\Gamma)$ is properly
equivalent to $\sigma\cdot\gamma(\Gamma)$.
}
\thmlab\FactorThrough
\proof
Let $b_L,b_R$ be the restrictions of $b:\Gamma\ra b\Gamma$
to $\Gamma_L,\Gamma_R$ respectively.
Then, as before, we have isometries $b_L:\Gamma_L\ra b\Gamma_L=(b\Gamma)_L$,
$b_R:\Gamma_R\ra b\Gamma_R=(b\Gamma)_R$, and they induce
isometries of discriminant groups
$$\eqalign{
&b_R^*:\Gamma_R^*/\Gamma_R\ra b\Gamma_R^*/b\Gamma_R\cr
&b_L^*:\Gamma_L^*/\Gamma_L\ra b\Gamma_L^*/b\Gamma_L.
}$$
By Lemma \Overlattice, we have
$$
\Gamma/(\Gamma_L\oplus\Gamma_R)
=\{\varphi(a)\oplus a|a\in\Gamma_R(-1)^*/\Gamma_R(-1)\}.
$$
Applying $b$ to this, we get
$$\eqalign{
b\Gamma/(b\Gamma_L\oplus b\Gamma_R)
&=\{b_L^*\circ\varphi(a)\oplus b_R^*(a)|a\in\Gamma_R(-1)^*/\Gamma_R(-1)\}\cr
&=\{b_L^*\circ\varphi\circ (b_R^*)^{-1}(c)\oplus c|c\in b\Gamma_R(-1)^*/b\Gamma_R(-1)\}.
}$$
It follows that $b_L^*\circ\varphi\circ (b_R^*)^{-1}:A_{b\Gamma_R(-1)}\ra A_{b\Gamma_L}$ is the
unique isomorphism determined by
the overlattice $b\Gamma$ of $b\Gamma_L\oplus b\Gamma_R$.
As a result, the triples
$$
\gamma(\Gamma)=(\Gamma_L,\Gamma_R(-1),\varphi),~~~
\gamma(b\Gamma)=(b\Gamma_L,b\Gamma_R(-1),b_L^*\circ\varphi\circ (b_R^*)^{-1})
$$
are properly equivalent. This completes the proof of assertion i.

We now consider the involution $\sigma:\cR\ra\cR$.
For convenience, we put
$$
\gamma(\Gamma)=(\Gamma_L,\Gamma_R,\varphi).
$$
Recall that
$$
(\sigma\Gamma)_L:=(\sigma\Gamma)\cap V=\sigma(\Gamma\cap V^\perp)=\sigma\Gamma_R,~~~
(\sigma\Gamma)_R:=(\sigma\Gamma)\cap V^\perp=\sigma\Gamma_L.
$$
Thus $\sigma$ induces a $\Z$-module isomorphism
$$
\Gamma/(\Gamma_L\oplus\Gamma_R)\br\sim\over\ra
\sigma\Gamma/((\sigma\Gamma)_R\oplus(\sigma\Gamma)_L),
~~(x,y)~mod~\Gamma_L\oplus\Gamma_R\mapsto
(\sigma x,\sigma^{-1}y)~mod~(\sigma\Gamma)_R\oplus(\sigma\Gamma)_L,
$$
and lattice isometries defined by
$$\eqalign{
&\iota_L:\Gamma_R(-1)\ra(\sigma\Gamma)_L,~~(0,y)\mapsto(\sigma^{-1}y,0)\cr
&\iota_R:\Gamma_L\ra(\sigma\Gamma)_R(-1),~~(x,0)\mapsto(0,\sigma x).
}$$
In turn, these isometries induce isometries of discriminant groups
$$\eqalign{
&\iota_L^*:\Gamma_R(-1)^*/\Gamma_R(-1)\ra{(\sigma\Gamma)_L}^*/(\sigma\Gamma)_L\cr
&\iota_R^*:\Gamma_L^*/\Gamma_L\ra(\sigma\Gamma)_R(-1)^*/(\sigma\Gamma)_R(-1).
}$$

By Lemma \Overlattice, we have
$$
\Gamma/(\Gamma_L\oplus\Gamma_R)=\{\varphi(a)\oplus a|a\in \Gamma_R(-1)^*/\Gamma_R(-1)\}.
$$
Applying $\sigma$ to this, we get
$$
\sigma\Gamma/((\sigma\Gamma)_R\oplus(\sigma\Gamma)_L)=\{\iota_L^*(a)\oplus\iota_R^*\circ\varphi(a)|
a\in \Gamma_R(-1)^*/\Gamma_R(-1)\}.
$$
Reparameterizing this set by setting $a'=\iota_R^*\circ\varphi(a)\in(\sigma\Gamma)_R(-1)^*/(\sigma\Gamma)_R(-1)$,
we get
$$\eqalign{
\sigma\Gamma/((\sigma\Gamma)_R\oplus(\sigma\Gamma)_L)&=\{\varphi'(a')\oplus a'|
a'\in (\sigma\Gamma)_R(-1)^*/(\sigma\Gamma)_R(-1)\}\cr
\varphi'&:=\iota_L^*\circ\varphi^{-1}\circ(\iota_R^*)^{-1}.
}$$
This shows that the triple $((\sigma\Gamma)_L,(\sigma\Gamma)_R(-1),\varphi')$
coincides with $\gamma(\sigma\Gamma)$, and
that it is properly equivalent to
$(\Gamma_R(-1),\Gamma_L,\varphi^{-1})=:\sigma\cdot\gamma(\Gamma)$.
This completes the proof of assertion ii.
$\Box$

\lemma{Let $\Gamma,\Gamma'\in\cR$.
\item{i.} If $\gamma(\Gamma')$ is properly equivalent to $\gamma(\Gamma)$,
then $\Gamma'=b\Gamma$ for some $b\in O(V)\times O(V^\perp)$.
\item{ii.} If $\gamma(\Gamma')$ is properly equivalent to $\sigma\cdot\gamma(\Gamma)$,
then $\Gamma'=b\sigma\Gamma$ for some $b\in O(V)\times O(V^\perp)$.
}\thmlab\Converse
\proof
Let's prove ii. assuming i. first. By the preceding lemma,
$\gamma(\sigma\Gamma)$ is properly equivalent to $\sigma\cdot\gamma(\Gamma)$.
Thus by hypothesis of ii., $\gamma(\Gamma')$ is
properly equivalent to $\gamma(\sigma\Gamma)$. It follows from i. that
$\Gamma'=b\sigma\Gamma$ for some $b\in O(V)\times O(V^\perp)$.

We prove i. now. Write
$$
\gamma(\Gamma)=(\Gamma_L,\Gamma_R,\varphi),~~~
\gamma(\Gamma')=(\Gamma_L',\Gamma_R',\varphi').
$$
By hypothesis of i., we have isometries
$$
b_L:\Gamma_L\ra\Gamma_L',~~~
b_R:\Gamma_R\ra\Gamma_R'
$$
such that
$$
\varphi'=b_L^*\circ\varphi\circ (b_R^*)^{-1}.
$$
There is a unique linear extension $b:\R^{2n}\ra\R^{2n}$, i.e.
$b|_{\Gamma_L}=b_L$ and $b|_{\Gamma_R}=b_R$.

Since $b$ restricted to the rank $2n$ lattice $\Gamma_L\oplus\Gamma_R\subset\R^{n,n}$
is an isometry, it follows that $b\in O(\R^{n,n})$.
Since $\Gamma_L\oplus\Gamma_R\subset\Gamma$, it follows that
$b(\Gamma_L\oplus\Gamma_R)=\Gamma_L'\oplus\Gamma_R'\subset b\Gamma$,
i.e. $b\Gamma$ is an overlattice of $\Gamma_L'\oplus\Gamma_R'$.
By Lemma \Overlattice, it determines a unique isometry
$\varphi_{b\Gamma}:A_{\Gamma_R'(-1)}\ra A_{\Gamma_L'}$.
By the same calculation as in the proof of Lemma \FactorThrough i.,
we see that
$$
\varphi_{b\Gamma}=b_L^*\circ\varphi\circ(b_R^*)^{-1}.
$$
Hence $\varphi_{b\Gamma}=\varphi'$. Thus $b\Gamma$ must
coincide with the unique overlattice corresponding to $\varphi'$, i.e.
$$
b\Gamma=\Gamma'.
$$

Finally, we have
$$
\Gamma'\cap V=\Gamma_L'=b\Gamma_L=b(\Gamma\cap V)
=b\Gamma\cap bV=\Gamma'\cap bV.
$$
This shows that both $n$ dimensional spaces $V$, $bV$
contains the rank $n$ lattice $\Gamma_L'$, hence must be identical: $bV=V$.
Likewise $bV^\perp=V^\perp$. It follows that $b\in O(V)\times O(V^\perp)$.
This completes the proof of assertion i. $\Box$

\theorem{The map $\cR\br\gamma\over\lra\cT$ descends to
injections
$$
\ol\cR\br\gamma\over\lra\ol\cT,~~~\ol\cR/\sigma\br\gamma\over\lra\ol\cT/\sigma.
$$
The image $\gamma(\ol\cR)$ consists of classes of triples $(X,Y,\varphi)$
where $X$ have the same signature as $V$.}
\thmlab\gammaTheorem
\proof
Lemma \FactorThrough~ says that
$\gamma$ descends to $\ol\cR\ra\ol\cT$
and $\ol\cR/\sigma\ra\ol\cT/\sigma$. Lemma \Converse~
says that the induced maps are injective.

Given a $V$-rational Narain lattice $\Gamma$, then
$\Gamma_L,\Gamma_R(-1)$ both have the same signature as $V$.
So, the triple $\gamma(\Gamma)$ has the same signature as $V$.
Conversely, if $(X,Y,\varphi)\in\cT$ has the signature of $V$, then
$X,Y(-1)$ can be respectively realized inside the vector
spaces $V,V^\perp$.
The isometry $\varphi:A_Y\ra A_X$
corresponds to a unique even unimodular overlattice $\Gamma\supset
X\oplus Y(-1)$ in $V\oplus V^\perp=\R^{n,n}$.
It is immediate that $\Gamma\in\cR$ and
$\gamma(\Gamma)=(X,Y,\varphi)$. This proves the second assertion.
$\Box$

\subsec{Isomorphisms of discriminants vs. triples}

Given lattices $X,Y$, let $Isom(A_X,A_Y)$ denote
the set of isometries of their discriminants.
Since an isometry of $X$ induces an isometry of $A_X$,
it follows that $O(X)$ acts from the left
on $Isom(A_X,A_Y)$. Likewise $O(Y)$ acts on
it from the right. The orbit space of this $O(X)\times O(Y)$ action
is the double quotient
$$
(*)~~~O(X)\backslash Isom(A_X,A_Y)/O(Y).
$$

\lemma{If the lattices $X,Y$ are not isomorphic, then the double quotient (*)
naturally parameterizes the set of improper equivalence classes of
triples $[X,Y,\varphi]\in\ol\cT/\sigma$.}

See Proposition 4.12 in [HLOYII]. The result there is
stated for positive definite binary quadratic form, but
the proof is valid in general.

\lemma{If the lattices $X,Y$ are isomorphic, primitive, and rank 2, then
the double quotient (*) naturally parameterizes the set of improper equivalence classes of
triples $[X,Y,\varphi]\in\ol\cT/\sigma$.}

See Proposition 4.12 in [HLOYII] for details.

\newsec{Primitive Sublattices}
\seclab\PiSection

\definition{
We denote by
$\Pi$ the set of rank $n$ primitive sublattices of the even unimodular
lattice $U^n$. Note that a primitive sublattice
need not be primitive as an abstract lattice.
We regard $O(U^n)$ as a transformation group of the set $\Pi$, and
write $\ol\Pi=\Pi/O(U^n)$.}

 In this section we will develop a correspondence between primitive
sublattices of $U^n$ and triples, which is parallel
to the correspondence we obtained, in the last section,
between $V$-rational Narain lattices and triples.

\bu {\it Anti-isometries.}
Let $\xi:U^n\ra U^n$ be an anti-isometry,
i.e. an automorphism of abelian groups such that
$$
\bra \xi(x),\xi(x)\ket=-\bra x,x\ket~~~\forall x\in U^n.
$$
For example, declaring $\xi:e_i\mapsto-e_i$
and $\xi:f_i\mapsto f_i$, defines an anti-isometry.
Obviously if $f$ is any isometry of $U^n$,
then $f\circ\xi$ and $\xi\circ f$ are both anti-isometries,
This shows that the set of all anti-isometries is
$\xi\cdot O(U^n)=O(U^n)\cdot\xi$,
and is independent of the choice of $\xi$.
Likewise if $\xi'$ is any other anti-isometry,
then $\xi'\circ\xi\in O(U^n)$, hence $\xi$ always behave
like an involution modulo $O(U^n)$.
It follows immediately that
$$
(*)~~~~\xi\cdot O(U^n)\cdot\xi'=O(U^n).
$$
Hence the group generated by anti-isometries
has index 2 over the group of isometries.
It is also straightforward to show that for each $M\in\Pi$, we have
$$
(\xi(M))^\perp=\xi(M^\perp).
$$

\lemma{Let $\xi$ be an anti-isometry of $U^n$. Then as
a group of transformations of $U^n$, the group
$\bra O(U^n),\xi\ket$ is independent of the choice of $\xi$.
Moreover, there is an induced action of this group on
the set $\Pi$, where $^\xi:M\mapsto M^\xi:=\xi(M)^\perp$.}
\thmlab\xiLemma
\proof
The first assertion follows immediately from the identity (*).
The group $O(U^n)$ acts on $\Pi$
by left translation.
To see that the group $\bra O(U^n),\xi\ket$ acts, it
suffices to show that for any $f\in O(U^n)$, and $M\in\Pi$, we have
the identities
$$
M^{f\circ\xi}=f(M^\xi),~~~~M^{\xi\circ f}=(fM)^\xi.
$$
But they readily follow from the definition of $M^\xi$. $\Box$

\corollary{
The set $\ol\Pi/\xi$ is independent
of the choice of $\xi$.}

From now on, we will fix an anti-isometry
$\xi$ of $U^n$ once and for all.

\lemma{Let $M\in\Pi$, and $f\in O(U^n)$. Then $fM$ and $M^\perp(-1)$
are isogenus to $M$, hence they have the same
signature, discriminant, and primitivity as $M$. The lattice $M^\xi\in\Pi$
has the same signature, discriminant, and primitivity as $M$.}
\thmlab\SDPLemma
\proof
Since $f$ is an isometry, $fM$ is isomorphic to $M$.
Now consider $M^\perp(-1)$. If $M$ has signature $(n-p,p)$,
then $M^\perp$ has signature $(p,n-p)$, since $U^n$
has signature $(n,n)$. It follows that $M^\perp(-1)$
has the same signature as $M$.
For discriminant, we have
$$
det~M^\perp(-1)=(-1)^n det~M^\perp=(-1)^n (-1)^n det~M=det~M.
$$
Since $U^n\supset M\oplus M^\perp(-1)$
is an even unimodular overlattice, it determines an isometry of
discriminant groups $\varphi(M):A_{M^\perp(-1)}\ra A_M$, by Lemma \Overlattice.
It follows that $M^\perp(-1)$ is isogenus to $M$, hence
has the same primitivity as $M$ by Lemma \BothPrimitive.

Next we consider $M^\xi$. That $M^\xi$ and $M$ have the same signature is
shown as in the case of $M^\perp(-1)$, and so we omit the details.
Since $\xi$ is an anti-isometry, it follows that $\xi(M)\in\Pi$
has the same primitivity as $M$.
Since $\xi(M)^\perp(-1)$ is isogenus to $\xi(M)$,
it follows that $M^\xi(-1)=\xi(M)^\perp(-1)$, hence $M^\xi$,
has the same primitivity as $\xi(M)$.  $\Box$

\definition{Define the map
$$
\theta:\Pi\ra\cT,
~~~~M\mapsto(M,M^\perp(-1),\varphi(M)).
$$
Here $\varphi=\varphi(M)$ is the isometry of discriminant groups
$A_{M^\perp(-1)}\ra A_M$ determined by the even unimodular
overlattice $U^n\supset M\oplus M^\perp$.}

Now comes the parallels of Lemmas \FactorThrough~
and \Converse.

\lemma{
\item{i.} If $M\in\Pi$ and $f\in O(U^n)$, then
$\theta(fM)$ is properly equivalent to $\theta(M)$.
\item{ii.} If $M\in\Pi$, then $\theta(M^\xi)$ is properly
equivalent to $\sigma\cdot\theta(M)$.
}
\thmlab\thetaFactorThrough
\proof
The proof is word for word the same as the proof of
Lemma \FactorThrough. The only slight difference is in ii.
Here it suffices to show that the triples
$$
\sigma\cdot\theta(M)=(M^\perp(-1),M,\varphi(M)^{-1}),~~~~
\theta(M^\xi)=(M^\xi,M^\xi~^\perp(-1),\varphi(M^\xi))
$$
are equivalent via the isometries
$$
\xi|_{M^\perp}:M^\perp(-1)\ra M^\xi,~~~~
\xi|_M:M\ra M^\xi~^\perp(-1).
$$
In other words, we need to show that
$$
\xi|_{M^\perp}^*\circ\varphi(M)^{-1}\circ \xi|_M^*~^{-1}=\varphi(M^\xi).
$$
This is a computation very similar to that in Lemma \FactorThrough.
We will not repeat it here. $\Box$

\lemma{Let $M,M'\in\Pi$.
\item{i.} If $\theta(M')$ is properly equivalent to $\theta(M)$,
then $M'=fM$ for some $f\in O(U^n)$.
\item{ii.} If $\theta(M')$ is properly equivalent to $\sigma\cdot\theta(M)$,
then $M'=f M^\xi$ for some $f\in O(U^n)$.
}\thmlab\thetaConverse
\proof
The proof is word for word the same as the proof of
Lemma \Converse, with the following minor change in the last part of i.
Being given a proper equivalence $\theta(M')\sim\theta(M)$ of triples
means that we are given isometries
$$
g:M\ra M',~~~~g_\perp:M^\perp\ra M'^\perp
$$
such that
$$
\varphi(M')=g^*\circ\varphi(M)\circ g_\perp^*~^{-1}.
$$
This implies that there is a unique isometry $f\in O(U^n)$ extending
$g$ and $g_\perp$. $\Box$

\theorem{The map $\Pi\br\theta\over\lra\cT$ descends to
bijections
$$
\ol\Pi\br\theta\over\lra\ol\cT,~~~\ol\Pi/\xi\br\theta\over\lra\ol\cT/\sigma.
$$
The maps preserves discriminant and primitivity.
}\thmlab\thetaTheorem
\proof
Lemma \thetaFactorThrough~ says that
$\theta$ descends to $\ol\Pi\ra\ol\cT$
and $\ol\Pi/\xi\ra\ol\cT/\sigma$. Lemma \thetaConverse~
says that the induced maps are injective.
To prove that the maps are surjective,
let $(X,Y,\phi)\in\cT$. This determines
a unique even unimodular overlattice $\Gamma\supset X\oplus Y(-1)$
of signature $(n,n)$. By Milnor's theorem,
there is an isometry $f:\Gamma\ra\Gamma^e$. This induces
a triple $\theta(f(X))=(f(X),f(Y)=f(X)^\perp(-1),\varphi(f(X)))$
which is properly equivalent to $(X,Y,\phi)$ via $f$.
This proves the asserted surjectivity.

Finally, by Lemma \SDPLemma, the map $\theta$ preserves discriminant and
(abstract lattice) primitivity. Since the equivalence relations defined
on $\Pi$ and on $\cT$ are compatible with primitivity and fixing discriminant,
the maps induced by $\theta$ must preserves these properties.
$\Box$

\newsec{Coprime Pairs and Gauss Product}

\subsec{Primitive sublattices and concordant pairs}

\bu {\it Notations.}
On the space $M^{2,2}$ of $2\times 2$ real matrices, define
the involutions
$$\eqalign{
&^t:\left[\matrix{\alpha&\beta\cr\gamma&\delta}\right]\mapsto
\left[\matrix{\alpha&\beta\cr\gamma&\delta}\right]^t:=
\left[\matrix{\alpha&\gamma\cr\beta&\delta}\right]\cr
&^A:\left[\matrix{\alpha&\beta\cr\gamma&\delta}\right]\mapsto
\left[\matrix{\alpha&\beta\cr\gamma&\delta}\right]^A:=
\left[\matrix{\delta&-\beta\cr-\gamma&\alpha}\right]\cr
&^\vee:\left[\matrix{\alpha&\beta\cr\gamma&\delta}\right]\mapsto
\left[\matrix{\alpha&-\beta\cr-\gamma&\delta}\right]\cr
&\sigma:\left[\matrix{\alpha&\beta\cr\gamma&\delta}\right]\mapsto
diag(1,-1)\left[\matrix{\alpha&\beta\cr\gamma&\delta}\right]
=\left[\matrix{\alpha&\beta\cr-\gamma&-\delta}\right].
}$$
Note that $^t$ and $^A$ are both anti-involution of the
matrix algebra $M^{2,2}$, i.e.
$$
(XY)^A=Y^AX^A,~~~(XY)^t=Y^tX^t.
$$
But $^\vee$ is an isometric involution of algebra,
i.e. $(XY)^\vee=X^\vee Y^\vee$, because
it is given by a conjugation:
$$
X^\vee=diag(1,-1)X diag(1,-1).
$$
All three are obviously isometries of the quadratic form $2det$ on $M^{2,2}$.
They are also pairwise commuting. The map $\sigma$ is
an anti-isometry:
$$
det~\sigma X=-det~X.
$$
We will use the fact that $M^{2,2}$ is isometric to $\R^{2,2}$.
See Appendix C in [HLOYII].

Recall that an element
$(g_1,g_2)\in P(SL(2,\R)^2)$ acts on $M^{2,2}$
by isometry via left and right multiplications
$X\mapsto g_1Xg_2^{-1}$.
Recall also that the subgroup of $P(SL(2,\R)^2)$
which stabilizes the lattice
$\Gamma^e\equiv\left[\matrix{\Z&\Z\cr\Z&\Z}\right]\subset M^{2,2}$
is $P(SL(2,\Z)^2)$. See Appendix C in [HLOYII] for details.

\lemma{The group $O(\Gamma^e)$ is generated by $P(SL(2,\Z)^2)$, $^A$, and $^\vee$.}
\thmlab\OGamma
\proof
As shown in Appendix C in [HLOYII], $O(\Gamma^e)$ has the shape
$\coprod_{g\in\cZ}g\cdot P(SL(2,\Z)^2)$ where $\cZ\subset O(\Gamma^e)$
is 4-element subgroup with 2 generators, such that
\item{i.} one generator is orientation reversing;
\item{ii.} one generator is orientation preserving,
but reversing a positive 2-plane orientation.

\noindent Note that $P(SL(2,\Z)^2)$ is orientation preserving,
and preserving a positive 2-plane orientation as well.
Then it is easy to show that $O(\Gamma^e)$
is generated by $P(SL(2,\Z)^2)$ plus {\it any two elements in $O(\Gamma^e)$} with
properties i.--ii. We verify that $^A\in O(\Gamma^e)$ has property i.,
and that $^\vee\in O(\Gamma^e)$ has property ii. This completes the proof. $\Box$


All quadratic forms are assumed {\it even non-degenerate}.
We always identify a quadratic form $P=[a,b,c]$ with its
matrix $\left[\matrix{2a&b\cr b&2c}\right]$, and write $disc~P=b^2-4ac$.
Recall that $\iota[a,b,c]:=[a,-b,c]$. If it is primitive,
then it represents
the inverse class of the class of $[a,b,c]$ in the class group.
Note that under the identification here, we have
$$
\iota[a,b,c]\equiv\left[\matrix{2a&b\cr b&2c}\right]^\vee.
$$
All three operations $^t,^A,^\vee$ operates on symmetric
matrices, hence on quadratic forms. We also write $-[a,b,c]=[-a,-b,-c]$.
If $P=[a,b,c]$, $Q=[a',b',c']$
are quadratic forms, we sometimes write
$$
gcd(P):=gcd(a,b,c),~~~~gcd(P,Q):=gcd(a,b,c,a',b',c').
$$

\definition{Let $P=[a,b,c],Q=[a',b',c']$ be forms
with the same discriminant, but not necessarily primitive.
Put $\delta:=gcd(a,b,c)$, $\delta':=gcd(a',b',c')$, $\lambda:=gcd(\delta,\delta')$.
We say that $P,Q$ are coprime if $\lambda=1$.
We say that $P,Q$ are concordant
if $aa'\neq0$, $gcd(a,a')=\lambda$, and $b=b'$. In this case,
we define
$$
P*Q:=[aa',b,{c\over a'}].
$$
(cf. Appendix A [HLOYII].)
If $d|\delta$, we also write ${1\over d}P=[{a\over d},{b\over d},{c\over d}]$.
We denote by $\cP$ the set of coprime pairs $(P,Q)$ of the same discriminant,
and by $\ol\cP$ the set of coprime pairs $(\bar P,\bar Q)$
of $GL(2,\Z)$ equivalence classes $\bar P,\bar Q$ of quadratic forms.
Let $\sigma:\cP\ra\cP$, $(P,Q)\mapsto(Q,P)$.}
\thmlab\CoprimePairs

Note that $P,Q$ being concordant implies that $ac=a'c'$.
Since $gcd(a,a')=\lambda$, it follows that ${a\over\lambda},{a'\over\lambda}$
are coprime. Thus we have ${a'\over\lambda}|{c\over\lambda}$
and ${a\over\lambda}|{c'\over\lambda}$. In particular we have ${c\over a'}\in\Z$.

\lemma{
Any pair of quadratic forms of the same discriminant
can be $SL(2,\Z)$ transformed to a concordant pair.}
\thmlab\ConcordantLemma

A proof can be found in Appendix A [HLOYII].

%
%

\subsec{The map $\Lambda$}

Recall that
$\Pi$ denotes the set of rank 2 primitive sublattices of
lattice $\Gamma^e\equiv\left[\matrix{\Z&\Z\cr\Z&\Z}\right]\subset M^{2,2}$ equipped with
the quadratic form $2det~X$.
Fix an anti-isometry $\xi:\Gamma^e\ra\Gamma^e$ once and for all.
We denote $\ol\Pi:=\Pi/O(\Gamma^e)$.

Let $P,Q$ be any (non-degenerate, as always)
binary quadratic forms. Introduce the notation
$$
\Lambda(P,Q):=\{X\in\Gamma^e|X^tP=QX^A\}.
$$
A priori, $\Lambda(P,Q)$
is just a primitive abelian subgroup of $\Gamma^e$.
But we will show shortly that it
is either zero or a rank 2 sublattice. In fact,
it has rank 2 iff $P,Q$ have the same discriminant.
The motivation for the defining equation for $\Lambda(P,Q)$ comes from the
fact that two quadratic forms $P,Q$ are $SL(2,\Z)$ equivalent iff there exists $g\in SL(2,\Z)$
such that $g^tPg=Q$, i.e. $g\in\Lambda(P,Q)$ a vector of length 2.

\lemma{For any quadratic forms $P,Q$,
the abelian group $\Lambda(P,Q)$ is nonzero iff $disc~P=disc~Q$.
In this case, $\Lambda(P,Q)$ has rank 2.}
\thmlab\RankTwo
\proof
Given $P=[a,b,c]$, $Q=[a',b',c']$, the equation $X^tP=QX^A$ becomes
$$
M\vec X=0,~~~M=\left[\matrix{
2 a& 0& b + b'& -2 a'\cr
b - b'& 2 a'& 2 c& 0\cr
0& 2 a& 2 c'& b - b'\cr
-2 c'& b + b'& 0& 2 c}\right].
$$
From this we find
$$
det~M=disc~Q-disc~P.
$$
This proves the first assertion.

We can now compute the $3\times 3$ minors of
$M$. Under the condition $det~M=0$, it is easy (with some help
from Mathematica) to check that all such minors are zero. This shows
that $M$ has rank at most 2. To see that it can
be no less than 2, observe that for $b=b'$, we have $rk~M=2$.
This shows that when $det~M=0$, then $rk~M$ is always 2,
implying that $rk~\Lambda(P,Q)=2$ as well.
$\Box$

\lemma{Let $P,Q$ be any quadratic forms. We have
\item{i.} $\Lambda(P,Q)^t=\Lambda(Q^A,P^A)$.
\item{ii.} $\Lambda(P,Q)^A=\Lambda(Q,P)$.
\item{iii.} $\Lambda(P,Q)^\vee=\Lambda(P^\vee,Q^\vee)$.
\item{iv.} $\Lambda(g^t Pg,Q)= g^{-1}\cdot\Lambda(P,Q)$ for $g\in SL(2,\Z)$.
\item{v.} $\Lambda(P,g^tQg)=\Lambda(P,Q)\cdot g$ for $g\in SL(2,\Z)$.
\item{vi.} $\Lambda(P^\vee,Q)= diag(1,-1)\cdot\Lambda(-P,Q)$.
\item{vii.} $\Lambda(P,Q^\vee)=\Lambda(P,-Q)\cdot diag(1,-1)$.
\item{viii.} $\Lambda(P,-Q)=\Lambda(-P,Q)$.
\item{ix.} $\Lambda(kP,kQ)=\Lambda(P,Q)$ for nonzero integer $k$.
\item{x.} $P\br SL(2,\Z)\over\sim Q$ iff $\exists X\in\Lambda(P,Q),~\bra X,X\ket=2det~X=2$.
\item{xi.} $P\br SL(2,\Z)\over\sim -Q^\vee$ iff $\exists X\in\Lambda(P,Q),~\bra X,X\ket=2det~X=-2$.
}\thmlab\InvarianceLemma
\proof
Assertions i.-x. are straightforward. To prove xi., notice
that by x., if $P\sim-Q^\vee$, then $\exists Y\in\Lambda(P,-Q^\vee)$ with $det~Y=2$.
But by vii., it follows that $X=Y\cdot diag(1,-1)\in\Lambda(P,Q)$ has $det~X=-2$.
The converse is similar. $\Box$

\lemma{For any quadratic forms $P,Q$ of discriminant $D$, the restriction of
$det$ to $\Lambda(P,Q)$ is nondegenerate. Moreover, $\Lambda(P,Q)$
is indefinite iff $P,Q$ are both indefinite.}
\proof
Since nondegeneracy is invariant
under isometries, we are free to replace $P,Q$ by their $SL(2,\Z)$
transforms, thanks to Lemma \InvarianceLemma iv.-v. So let's
assume that $P=[a,b,c]$, $Q=[a',b',c']$,
are concordant. Then the coefficient matrix $M$ in Lemma \RankTwo~
simplifies, and finding a $\Q$-basis of $\Q~\Lambda(P,Q)$ is easy.
The result is
$$
\Q~\Lambda(P,Q)=\Q\left\{\right.
\left[\matrix{{a'\over a}&0\cr 0&1}\right],
\left[\matrix{-{b\over a}&-{c\over a'}\cr 1&0}\right]\left.\right\}.
$$
Computing the discriminant of the quadratic form
on $\Q~\Lambda(P,Q)$ using this $\Q$-basis, we get
the result $D\over a^2$, where $D=b^2-4ac$.
This shows that $\Lambda(P,Q)$ is nondegenerate,
and that $\Lambda(P,Q)$ is indefinite iff both $P,Q$ are. $\Box$

\definition{Recall that $\cP$ denotes
the set of coprime pairs of binary quadratic forms
having the same discriminant. Define a map
$$
\Lambda:\cP\ra\Pi,~~~(P,Q)\mapsto\Lambda(P,Q).
$$}

\lemma{Let $P,Q$ be any quadratic forms of the same discriminant $D$.
If $P,Q$ are both positive definite or both negative definite,
then $\Lambda(P,Q)$ is positive definite.  If $P$ is positive definite and
$Q$ is negative definite (or the other way around),
then $\Lambda(P,Q)$ is negative definite.}
\thmlab\SignatureLemma
\proof
By the preceding lemma, if $P,Q$ are definite,
then $\Lambda(P,Q)$ is also definite. So to
determine the sign, it suffices to compute
the length of a single nonzero vector $X\in\Lambda(P,Q)$.
We have
$$
X^tP=QX^A.
$$
Since $X^A=(det~X)X^{-1}$,
we get an equivalent equation
$$
X^tPX=(det~X)Q.
$$
Note that $X^t P X$ have the same signature as $P$.
So if $P,Q$ are both positive definite or both negative definite,
then $det~X$ must be positive, and hence $\Lambda(P,Q)$
is positive definite. If $P$ is positive definite and
$Q$ is negative definite, then $det~X$ is negative and $\Lambda(P,Q)$
is negative definite. $\Box$

\lemma{For any quadratic forms $P,Q$ of the same discriminant, we have
$$
\Lambda(P,Q)^\perp=\Lambda(P,-Q)=(~^t\circ\sigma\circ~^t)\cdot\Lambda(P,Q^\vee)
$$}
\thmlab\PerpLemma
\proof
The second equality follows from that
$$
(~^t\circ\sigma\circ~^t)Z=Z~diag(1,-1),~~\forall Z\in M^{2,2},
$$
and Lemma \InvarianceLemma vii. that $\Lambda(P,Q^\vee)\cdot diag(1,-1)=\Lambda(P,-Q)$.

Let's consider the first equality.
Since both $\Lambda(P,-Q),\Lambda(P,Q)$
are rank 2 primitive sublattices in $\Gamma^e$,
it suffices to show that if $X\in\Lambda(P,-Q)$
and $Y\in\Lambda(P,Q)$, then $\bra X,Y\ket=0$.
So suppose that
$$
X^tP=-QX^A,~~~~Y^tP=QY^A.
$$
Then we have
$$
(X+Y)^tP=-Q(X-Y)^A.
$$
Taking the determinant on both sides, and noting that $det~P=det~Q\neq0$, we see that
$det(X+Y)=det(X-Y)$. This shows that
$$
\bra X,Y\ket=\bra X,-Y\ket=-\bra X,Y\ket.
$$
Here we've used the fact that $X,Y$ are $2\times 2$.
This completes the proof. $\Box$

\lemma{Suppose that $P=[a,b,c]$, $Q^\vee=[a',b,c']$ are
concordant forms with $gcd(a,a')=1$. Then
$$
\Lambda(P,Q)=\Z\left\{\right.
\left[\matrix{a'&-b\cr 0&a}\right],
\left[\matrix{0&-{c\over a'}\cr 1&0}\right]\left.\right\}.
$$}
\thmlab\LambdaBasis
\proof
Name the two given generators $X_1,X_2$.
It is trivial to check that the $X_i$
solve the linear equation
$X^t P-QX^A=0$. So the $X_i$ form a
$\Q$-basis of $\Lambda(P,Q)$. To show that the $X_i$
form a $\Z$-basis, it suffices to show that
the vectors $X_i\in\Z^4$
generate a primitive sublattice of $\Z^4$. Since
$X_1=(a',-b,a,0)$,
$X_2=(0,-{c\over a'},0,1)$,
we can find two additional vectors $Y_1,Y_2\in\Z^4$ of the shape $(*,*,*,0)$
such that the four vectors $X_1,X_2,Y_1,Y_2$
form a unimodular matrix. This shows that the first two vectors
generates a primitive sublattice of $\Z^4$.
$\Box$
\corollary{Let $P=[a,b,c], Q^\vee=[a',b',c']$ be quadratic forms of discriminant $D$.
Put  $\delta:=gcd(P)$, $\delta':=gcd(Q)$, $\lambda:=gcd(\delta,\delta')$, as before.
If $P, Q^\vee$ are concordant, then
$$
\Lambda(P,Q)=\Z\left\{\right.
{1\over\lambda}\left[\matrix{a'&-b\cr 0&a}\right],
\left[\matrix{0&-{c\over a'}\cr 1&0}\right]\left.\right\}.
$$
In particular, it has discriminant ${1\over{\lambda}^2}(b^2-4ac)$.
More generally if $P,Q$ are arbitrary, then
$\Lambda(P,Q)$ has discriminant ${D\over{\lambda}^2}$.}
\thmlab\MostGeneral
\proof By definition, that $P,Q^\vee$ are concordant means
that $aa'\neq0$, $gcd(a,a')=\lambda$, and $b=b'$.
It follows that ${1\over\lambda}P,~{1\over\lambda}Q^\vee$
are also concordant, but with coprime and nonzero
leading coefficients ${a\over\lambda},{a'\over\lambda}$.
By Lemma \InvarianceLemma ix., we have
$$
\Lambda(P,Q)=\Lambda({1\over\lambda}P,{1\over\lambda}Q).
$$
Now applying Lemma \LambdaBasis~ to compute
the right hand side, we get our first assertion.
Computing the discriminant using the $\Z$-base we found
is straightforward, and we get ${1\over{\lambda}^2}(b^2-4ac)$.

Finally, if $P,Q^\vee$ are arbitrary forms of discriminant $D$,
we can replace them with their $SL(2,\Z)$ transforms without
changing the values of $\delta,\delta',\lambda$,
or the discriminant of $\Lambda(P,Q)$,
thanks to Lemma \InvarianceLemma iv.-v. We can choose the
$SL(2,\Z)$ tranforms to be a concordant pair,
by Lemma \ConcordantLemma. Now our third assertion follows
from the second assertion. $\Box$

\corollary{$\Lambda(P,Q)$ is primitive as an abstract lattice iff
${1\over\lambda}P,{1\over\lambda}Q$ are both primitive as quadratic forms.}
\thmlab\LambdaPrimitivity
\proof
Again, isometries preserve primitivity. Thus,
we may as well assume that $P=[a,b,c],Q^\vee=[a',b,c']$ are concordant,
thanks to Lemmas \InvarianceLemma iv.-v. and the existence of concordant forms.
So we can use the $\Z$-base found in the preceding corollary
to compute the quadratic form of $\Lambda(P,Q)$. Since
$$
\Lambda(P,Q)=\Lambda({1\over\lambda}P,{1\over\lambda}Q),
$$
we may as well assume that $\lambda=1$.
Then the resulting quadratic form is $[aa',b,{c\over a'}]$.
We put
$$
m:=gcd(aa',b,{c\over a'}).
$$

Suppose $P,Q$ are both primitive.
Then we have $m|aa'$, $m|b$, and $m|{c\over a'}$, hence $m|c$.
Since $[a,b,c]$ is assumed primitive i.e. $gcd(a,b,c)=1$, it follows that
$m|a'$. Since $a'c'=ac$, and $m|{c\over a'}$,
it follows that $m|{c'\over a}$, hence $m|c'$.
So we find that $m|a',b,c'$.
But $[a',b,c']$ is also assumed primitive. It follows
that $m=1$, and so the quadratic form of $\Lambda(P,Q)$ is primitive.

Conversely suppose that one of $P,Q$ is not primitive.
By Lemma \InvarianceLemma ii., we may assume that $P$ is not primitive.
So $\delta:=gcd(a,b,c)>1$. Since $P,Q$ are assumed concordant,
we have $gcd(a,a')=1$, hence $gcd(\delta,a')=1$. It follows that
$\delta|{c\over a'}$. But we also have $\delta|a,b$.
So we have $\delta|m=gcd(aa',b,{c\over a'})$.
This shows that $m>1$, hence
$\Lambda(P,Q)$ is not primitive. $\Box$

\theorem{If $P,Q^\vee$ are coprime concordant forms,
then the lattice $\Lambda(P,Q)$
equipped with the $\Z$-base as in Lemma \LambdaBasis,
coincides with the quadratic form $[aa',b,{c\over a'}]=P*Q^\vee$.
More generally, for any coprime forms $P,Q$,
the lattices $\Lambda(P,Q)$, $\Lambda(P,Q)^\perp(-1)$,
are respectively isomorphic to $P*Q^\vee$, $P*Q$ as lattices.}
\thmlab\LambdaGauss
\proof
The first assertion follows from a straightforward computation using
the given explicit $\Z$-base. That $\Lambda(P,Q)$ is isomorphic to
$P*Q^\vee$ follows from the first assertion and Lemma \InvarianceLemma iv.-v.

Note that the transformation
$\xi:={}^t\circ\sigma\circ{}^t:\Gamma^e\ra\Gamma^e$
is an anti-isometry. In particular  for
any sublattice $M\subset\Gamma^e$,
we have $\xi(M)(-1)\cong M$ as lattices.
By Lemma \PerpLemma, it follows that
$\Lambda(P,Q)^\perp(-1)$ is isomorphic to $\Lambda(P,Q^\vee)$.
The latter is isomorphic to $P*Q$, by the first assertion.
$\Box$

\corollary{When restricted to primitive forms of discriminant $D$,
the map $\Lambda$ descends to
$$
\Lambda:Cl_D\times Cl_D\ra Cl_D/\iota,~~~(P,Q)\mapsto[\Lambda(P,Q)]=[P*Q^\vee].
$$
Here $Cl_D$ is the group of proper equivalence classes of
primitive forms (cf. p337 [Ca]), and $\iota:C\mapsto C^{-1}$.}
\proof
By Lemma \InvarianceLemma iv.-v.,
the isomorphism class $[\Lambda(P,Q)]$ depends only on
the $SL(2,\Z)$ equivalence classes of $P,Q$, i.e.
$\Lambda$ is a class function.
We evaluate $\Lambda(P,Q)$ by chosing
$P,Q^\vee$ to be concordant. Now the
preceding corollary completes the proof. $\Box$

If we restrict ourselves to positive definite primitive forms,
then it can be shown that there is a natural lifting to
$$
\Lambda:Cl_D\times Cl_D\ra Cl_D,~~~(C_1,C_2)\mapsto(C_1*C_2^{-1}).
$$
The point is that because the Grassmannian
of {\it positive} 2-planes is contractible,
we can choose an orientation for every positive 2-plane
so that they are all compatible under deformation.
In particular, we can assign compatible orientations to all $\Lambda(P,Q)$,
so that these lattices become quadratic forms. There are obviously
two ways to do so. The observation here is that
one of them yields the map $(C_1,C_2)\mapsto(C_1*C_2^{-1})$,
and the other choice yields $(C_1,C_2)\mapsto(C_1^{-1}*C_2)$.


Define the following transformations on the set $\cP$:
$$\eqalign{
P(SL(2,\Z)^2)\ni(g_1,g_2)&:(P,Q)\mapsto (g_1^t~^{-1}Pg_1^{-1},g_2^t~^{-1}Qg_2^{-1})\cr
\sigma=\iota_1&:(P,Q)\mapsto(Q,P)\cr
\iota_2&:(P,Q)\mapsto(P^\vee,Q^\vee)\cr
\iota_3&:(P,Q)\mapsto(P,Q^\vee).
}$$
It is clear that each of these transformations preserves
discriminant, coprimeness, and primitivity (but not necessary signature).
We denote the group generated by these transformations by
$\G=\bra P(SL(2,\Z)^2),\iota_1,\iota_2,\iota_3\ket$.
It is straightforward to check that $\ol\cP=\cP/K$ (see Definition \CoprimePairs),
where $K\subset\G$ is
the subgroup generated by $P(SL(2,\Z)^2)$ and $\iota_2,\iota_3$.
Note that among the generators of the group $\G$ presented above,
$\iota_1$ is the only one which does not preserve the signatures
of a pair, since $(P,Q),(Q,P)\in\cP$ have different signatures
when $P$ is positive definite and $Q$ is negative definite.

By Lemma \OGamma, the group $O(\Gamma^e)$ is generated by $P(SL(2,\Z)^2),^A,^\vee$.
Together with the anti-isometry $\xi$, they
act on primitive sublattices $M\in\Pi$ by
$$\eqalign{
P(SL(2,\Z)^2)\ni (g_1,g_2)&:M\mapsto g_1Mg_2^{-1}\cr
^A&:M\mapsto M^A\cr
^\vee&:M\mapsto M^\vee\cr
^\xi&:M\mapsto M^\xi.
}$$
By Lemma \SDPLemma, each of these transformations preserves discriminant,
and (abstract lattice) primitivity, and signature.

\theorem{The map $\Lambda:\cP\ra\Pi$
preserves discriminant, primitivity, and is equivariant
with respect to the group action of $\G$ on $\cP$,
and $\bra O(\Gamma^e),\xi\ket$ on $\Pi$.
Thus $\Lambda$ descends to $\ol\cP/\sigma\ra\ol\Pi/\xi$.}
\thmlab\LambdaEquivariance
\proof
If $(P,Q)$ is a coprime pair of discriminant $D$,
then $\Lambda(P,Q)$ has discriminant $D$, by
Corollary \MostGeneral, and so $\Lambda$ preserves discriminant.
If $(P,Q)$ is a primitive pair of discriminant $D$,
then $\Lambda(P,Q)$ is primitive, by Lemma \LambdaPrimitivity,
and so $\Lambda$ preserves primitivity.

Moreover, by Lemmas \InvarianceLemma iv.-v., ii.-iii., the map $\Lambda$
is equivariant with respective to $P(SL(2,\Z)^2)$ and the involutions
$\sigma=\iota_1\leftrightarrow ^A$, $\iota_2\leftrightarrow ^\vee$.
Since $\ol\Pi/\xi=\Pi/\bra O(\Gamma^e),\xi\ket$ is
independent of the choice of the anti-isometry $\xi$, by Lemma \xiLemma,
we can choose $\xi$ to be
$$
\xi=~^t~\circ\sigma\circ~^t.
$$
Then Lemmas \PerpLemma~ and \xiLemma~ show that the map $\Lambda$ is
also equivariant with respect to the involutions
$\iota_3\leftrightarrow ^\xi$, i.e. $\Lambda(\iota_3(P,Q))=\Lambda(P,Q)^\xi$.
Thus we have shown that the map $\Lambda$
is equivariant with respect to the group action of $\G$ on $\cP$,
and $\bra O(\Gamma^e),\xi\ket$ on $\Pi$. $\Box$


\newsec{Surjectivity of $\Lambda$}

\lemma{Let $A\in GL(2,\Q)$ which is not a multiple of the identity. Let
$M^{2,2}_\Q=\left[\matrix{\Q&\Q\cr\Q&\Q}\right]$.
Then the linear map $L_A:M^{2,2}_\Q\ra M^{2,2}_\Q$, $P\mapsto A^tPA-(det~A)P$
has rank 2, and we have
$$
ker~L_A=\Q\left\{\right.
\left[\matrix{0&1\cr-1&0}\right]
\left[\matrix{-c&a-d\cr0&b}\right]\left.\right\},~~~~
A=\left[\matrix{a&b\cr c&d}\right].
$$}
\proof
This is straightforward linear algebra. $\Box$

\corollary{Under the same assumption as in the preceding lemma,
there is a unique, up to sign, primitive quadratic form in $ker~L_A$.}
\proof
It is obvious that the 2-dimensional space $ker~L_A$
contains nonsymmetric matrices. So the subspace of
symmetric matrices is at most 1-dimensional. But $ker~L_A$ is
closed under transpose. It follows that the
subspace of symmetric matrices
in $ker~L_A$ is exactly 1-dimensional.
Since this is over $\Q$, there is a unique, up to sign,
primitive vector
$$
\left[\matrix{x&y\cr y&z}\right]\in\left[\matrix{\Z&\Z\cr\Z&\Z}\right]
$$
in this subspace, i.e. $gcd(x,y,z)=1$. If both $x,z$ are even,
then $\pm \left[\matrix{x&y\cr y&z}\right]$ are the
quadratic forms we seek. If not, then
$\pm 2\left[\matrix{x&y\cr y&z}\right]$ are.
$\Box$


\corollary{Let $X_1,X_2\in\Gamma^e$, be any linearly independent
vectors with nonzero lengths. Then there is a unique, up to overall sign,
pair of coprime forms $P,Q$ of the same discriminant such that
$X_1,X_2\in\Lambda(P,Q)$.}
\proof
Given the $X_i$, we want to solve the equations
$$
X_1^tP-QX_1^A=0,~~~~X_2^tP-QX_2^A=0.
$$
Since $det~X_i\neq0$ by assumption, it follows that $X_i^A=(det~X_i)X_i^{-1}$.
Moreover any solution $P,Q$ necessarily have $(det~X_i)(det~P)=(det~Q)(det~X_i)$,
i.e. $P,Q$ have the same discriminant.
We also get equivalent equations
$$
(*)~~~~X_1^tPX_1-(det~X_1)Q=0,~~~~X_2^tPX_2-(det~X_2)Q=0.
$$
Eliminating $Q$, we get
$$
A^t P A-(det~A)P=0,~~~~A:=X_1X_2^{-1}.
$$
Note that $A\in GL(2,\Q)$, and that
$A$ is not a scalar multiple of the identity,
because $X_1,X_2$ are assumed linearly independent.
So there is a unique, up to sign,
primitive quadratic form $P$ solving the last equation.
Plugging $P$ back into either equation in (*)
and solve for $Q$, we get
$$
Q={1\over det~X_1}X_1^tPX_1={1\over det~X_2}X_2^tPX_2.
$$
It is symmetric and nondegenerate because $P$ is. Its entries are rational numbers.
So by a suitable scaling of $P$ by an integer,
we get an even integral quadratic form $Q$.

From the construction, it is clear that $P,Q$ can be made coprime.
Uniqueness is also clear. $\Box$

\corollary{Any rank 2 primitive sublattice $M\subset\Gamma^e$
can be realized as $\Lambda(P,Q)$ for a pair of coprime
quadratic forms $P,Q$ of the same discriminant. Moreover $P,Q$
are unique up to an overall sign.}
\proof
By assumption, we can find two linearly
independent vectors $X_1,X_2\in M$ with nonzero lengths (nonzero lengths
because $det$ restricted to $M$ is assumed nondegenerate).
By the preceding corollary, we can find $P,Q$
such that $X_1,X_2\in\Lambda(P,Q)$.
Since $M$ and $\Lambda(P,Q)$ are both primitive sublattices of $\Gamma^e$
containing $X_1,X_2$, it follows that $M=\Lambda(P,Q)$.

Now a priori, $P,Q$ may depend on the choice of $X_1,X_2$.
But if $\Lambda(P,Q)=\Lambda(P',Q')$ for two coprime pairs $(P,Q),(P',Q')$,
then $X_1,X_2\in\Lambda(P',Q')$ also determine $P',Q'$
in terms of $X_1,X_2$ up to an overall sign. It follows that
$(P',Q')=\pm(P,Q)$.
$\Box$

\theorem{The map $\Lambda:\cP\ra\Pi$
is unramified 2:1 surjective, with fiber $\Lambda^{-1}(\Lambda(P,Q))=\{\pm(P,Q)\}$.}
\thmlab\TwoOneTheorem
\proof
By Lemma \InvarianceLemma viii., we have
$\Lambda(-P,-Q)=\Lambda(P,Q)=:M$.
The preceding corollary says $\Lambda$ is surjective
with fiber $\Lambda^{-1}(M)=\{\pm(P,Q)\}$.
Since $(P,Q),(-P,-Q)$
are never equal in $\cP$,
it follows that the map $\Lambda$ is 2:1 everywhere.
$\Box$

\theorem{The map $\Lambda:\cP\ra\Pi$ descends to
a bijection
$$
\Lambda:\ol\cP/\bra\sigma,-\ket\ra\ol\Pi/\xi
$$
where $\sigma,-$ are the respective involutions
on pairs $\sigma:(P,Q)\mapsto(Q,P)$, $-:(P,Q)\mapsto(-P,-Q)$.
Moreover, $\Lambda$ preserves discriminant and primitivity.}
\thmlab\LambdaTheorem
\proof
By Theorem \LambdaEquivariance, $\Lambda:\cP\ra\Pi$ is equivariant with respect to
$\G$ acting on $\cP$ and $K:=\bra O(\Gamma^e),\xi\ket$ acting on $\Pi$.
By Theorem \TwoOneTheorem, this map descends to an equivariant {\it bijection}
$$
\Lambda:\cP/\bra-\ket\ra\Pi.
$$
(Note that the $\G$ action commutes with $-$.)
Now passing to the orbit spaces of $\G$ and $K$, we get
the asserted induced bijection.

That $\Lambda$ preserves discriminant follows from Theorem \LambdaGauss.
That $\Lambda$ preserves primitivity follows from Corollary \LambdaPrimitivity.
The same holds for the induced $\Lambda$,
since the equivalence relations on $\cP$ and $\Pi$ are
compatible with primitivity and fixing discriminant.
$\Box$

\newsec{Applications and Conclusions}

\subsec{Binary forms}

Throughout this section, we set $n=2$.

Fix a negative integer $D$, and let
$\cP_D^+\subset\cP$ be the set of coprime pairs of positive
definite forms in $\cP$ of discriminant $D$. Let $\Pi_D^+\subset\Pi$ be the positive
definite lattices in $\Pi$ of discriminant $D$.
Then Theorem \LambdaTheorem~ says that
\eqn\dumbI{
\ol\cP_D^+/\sigma \br\Lambda~1:1\over\longleftrightarrow\ol\Pi_D^+/\xi.
}
The correspondence \dumbI~ also preserves primitivity. Note that
this correspondence can be precisely described by the composition law $*$,
as in Theorem \LambdaGauss.

Let $\cT_D^+$ be the set of positive definite triples in $\cT$ of discriminant $D$.
Then Theorem \thetaTheorem~ says that
\eqn\dumbII{
\ol\Pi_D^+/\xi \br\theta~1:1\over\longleftrightarrow\ol\cT_D^+/\sigma.
}
This correspondence also preserves primitivity.

Now put $V=\R^{2,0}$, and let $\cR^+_D$ be the set of $V$-rational Narain lattices
in $\cR$ of discriminant $D$. Then Theorem \gammaTheorem~ says that
\eqn\dumbIII{
\ol\cR^+_D/\sigma \br\gamma~1:1\over\longleftrightarrow\ol\cT_D^+/\sigma.
}
Note that the same is true if $V$ is replaced by
any other positive definite two-plane in $\R^{2,2}$.
The three correspondences above together recover a result in [HLOYII]
(cf. Theorem 5.8 there).

Let's consider now the indefinite case, i.e. $D$ a positive integer,
and $V=\R^{1,1}=\{(*,0,*,0)\}\subset\R^{2,2}$.
Let $\cP_D^{+-}$, $\Pi_D^{+-}$, $\cT_D^{+-}$, $\cR^{+-}_D$,
be respectively the set of coprime pairs of indefinite forms in $\cP$,
the set of indefinite lattices in $\Pi$, the set of
indefinite triples in $\cT$, and finally, the set
of $V$-rational Narain lattices in $\cR$, all having discriminant $D$.
Then the correspondences analogous to \dumbII~ and \dumbIII~
hold. But \dumbI~ must be replaced by
\eqn\dumbIV{
\ol\cP_D^{+-}/\bra \sigma,-\ket \br\Lambda~1:1\over
\longleftrightarrow\ol\Pi_D^{+-}/\xi.
}
While in the definite case, the involution
$-:(P,Q)\mapsto(-P,-Q)$ always identifies two distinct
pairs in $\ol\cP$, it is no longer so in the indefinite case.
The reason is that a given indefinite form $P$
may or may not be $GL(2,\Z)$ equivalent to $-P$. Both possibilites
can occur.

For simplicity, let's consider just the primitive indefinite forms.
There are two kinds of $GL(2,\Z)$ classes. There are classes with
$\bar P=-\bar P$, and those with $\bar Q\neq-\bar Q$. Clearly the
latter kind comes in pairs. Thus the complete list
of pairwise distinct $GL(2,\Z)$ classes has the shape
$$
\bar P_1,...,\bar P_u,\bar Q_1,...,\bar Q_v,-\bar Q_1,...,-\bar Q_v.
$$
The numbers $u,v$ depend on $D$.

To describe the primitive objects in the set $\ol\cP_D^{+-}/\bra\sigma,-\ket$,
it suffices to list all pairs which are not congruent
modulo $\bra\sigma,-\ket$. We get
$$\eqalign{
&(\bar P_i,\bar P_j),~~~i\leq j\cr
&(\bar P_i,\bar Q_j),~~~\forall i,j\cr
&(\bar Q_i,\bar Q_j),~~~i\leq j\cr
&(\bar Q_i,-\bar Q_j),~~~\forall i,j.
}$$
This shows that
$$
\#primitive~objects~in~\ol\cP_D^{+-}/\bra\sigma,-\ket
=\half u(u+1)+uv+\half v(v+1)+v^2=\half(u+v)(u+v+1)+v^2.
$$
It follows that this is also the number of primitive objects in
each of $\ol\Pi_D^{+-}/\xi, \ol\cT^{+-}_D/\sigma, \ol\cR^{+-}_D/\sigma$.

It is known that for $D=p$ a prime number,
every quadratic form $P$ of discriminant $D$
is $GL(2,\Z)$ equivalent to $-P$. Thus, in this case,
we have $v=0$.

\subsec{Concluding remarks}

%

We now return to the general setting of {\it The Main Problem}
in the Introduction. Here $\R^{n,n}$ is replaced by
an arbitrary quadratic space $E$ of signature $(r,s)$ with $8|(r-s)$
(which can be taken to be $\R^{r,s}$ without loss of generality),
and $V\subset E$ an arbitrary non-degenerate subspace
(which can be taken to be $\R^{p,q}:=
\{(*^p,0^{r-p},*^q,0^{s-q})\}\subset\R^{r,s}$).
It is further assumed that $E$ is indefinite.
(Note that Milnor's theorem for uniqueness
of indefinite unimodular lattices is used in Theorem \thetaTheorem.)
Then Theorems \gammaTheorem~ and \thetaTheorem~
can be readily generalized to this case.
The notion of a triple must be weakened as follows: {\it A triple
$(X,Y,\varphi)$ consists of a pair of even lattices $X,Y$,
and an isometry $\varphi:A_Y\ra A_X$,
such that $X\oplus Y(-1)$ have signature $(r,s)$.}
The notions of $V$-equivalence on $V$-rational
Narain lattices and proper equivalence
on triples remain the same as before.
Thus $\cR$ is now the set of $V$-rational Narain lattices,
$\cT$ the class of triples, and $\ol\cR$, $\ol\cT$ the
respective sets of equivalence classes.
The set $\Pi$ now becomes the set of primitive sublattices
of a given abstract even unimodular lattice $\bf U$ of signature $(r,s)$,
and $\bar\Pi$ the set of $O({\bf U})$ equivalence classes in $\Pi$.
The lattice $\bf U$ now plays the role of the lattice $U^n$.

In this generality, however, we no longer
have an anti-isometry $\sigma$ of $E$ which
exchanges $V$ and $V^\perp$. Likewise, for $M\in\Pi$,
the lattice $M^\perp(-1)$ is no longer isogenus to $M$ in general.
In any case, the set of proper equivalence classes
of triples $(X,Y,\varphi)$ with $\varphi\in Isom(A_Y,A_X)$,
and $X,Y$ are lattices of fixed discriminant $D$
with $sign~X\oplus Y(-1)=(r,s)$,
is clearly finite. Therefore, the problem of counting
$V$-rational Narain lattices,
triples, and primitive sublattices of $\bf U$, still makes sense.
In fact, Theorems \gammaTheorem~ and \thetaTheorem~ generalize to

\theorem{The map $\cR{\br\gamma\over\ra}\cT$, $\Gamma\mapsto
(\Gamma\cap V,(\Gamma\cap V^\perp)(-1),\varphi_\Gamma)$,
descends to an injection $\ol\cR\br\gamma\over\ra\ol\cT$. The image
$\gamma(\ol\cR)$ consists of classes of triples $(X,Y,\varphi)$ where
$X$ has the same signature as $V$.}

\theorem{The map $\Pi{\br\theta\over\ra}\cT$,
$M\mapsto(M,M^\perp(-1),\varphi(M))$, descends to a bijection
$\ol\Pi\br\theta\over\ra\ol\cT$.}

The proofs of Theorems \gammaTheorem~ and \thetaTheorem~
are valid here verbatim. We give two applications of these theorems.

Let $E=\R^{1,1}$, and $V\subset E$ a 1
dimensional positive definite subspace. We want to count
$V$-rational Narain lattices $\Gamma$ with $\Gamma\cap V$
rank 1 of discriminant $D>0$. This amounts to counting
rank 1 primitive sublattices $M\subset U$ of discriminant $D$
modulo $O(U)$.
It is not hard to show that the answer is 1 if $D=2$, and
exactly $2^{\tau({D\over 2})-1}$ if $D>2$.
Here $\tau(k)$ is the number of distinct prime factors of $k$.

We now consider the general case $E=\R^{r,s}$, and $V=\R^{p,q}\subset E$.
As before, denote by $\cR_D$
the set of $V$-rational Narain lattices
$\Gamma$ such that $\Gamma\cap V$ has discriminant $D$,
and by $\ol\cR_D$ the set of $V$-equivalence classes in $\cR_D$.
Let $\cL_D$ denotes the class of abstract even lattices $L$
of discriminant $D$, and $\ol\cL_D$ the set of
isomorphism classes $[L]$ in $\cL_D$. Note that $\ol\cL_D$
decomposes into finite subsets (see [Ca]) $\ol\cL_D(p,q)$
consisting of isomorphism classes of lattices
of discriminant $D$ and signature $(p,q)$.

\proposition{If $r>p$, $s>q$,
and $dim~V\neq\half dim~E$, then $|\ol\cR_D|=|\ol\cL_D(p,q)|$.}
\proof
Note that a Narain lattice $\Gamma\subset E$ is $V$-rational
iff it is also $V^\perp$-rational. Thus it suffices
to consider the case $dim~V<\half dim~E$, since $dim~V+dim~V^\perp=dim~E$.
By the preceding theorems, it suffices to show that the map
$$
(*)~~~~\ol\Pi_D(p,q)\ra\ol\cL_D(p,q),~~~M~mod~O({\bf U})\mapsto[M],
$$
is a bijection, where $\ol\Pi_D(p,q)$ is the set consisting
of the objects of discriminant $D$ and signature $(p,q)$ in $\ol\Pi$.

Let $[M]=[M']$ for some primitive sublattices
$M,M'\in\Pi$ of $\bf U$ with discriminant $D$,
signature $(p,q)$, and let $i:M\hra{\bf U}$,
$i':M'\hra{\bf U}$ be their respective inclusions.
Then we have a lattice isomorphism $M{\br g\over\ra}M'$,
and two primitive embeddings $i,i'\circ g:M\ra{\bf U}$.
Since $p+q=rank~M=dim~V\leq\half dim~E-1=\half(r+s)-1$, it follows that
$$
r+s-p-q\geq 2+rank~M\geq 2+l(A_M)
$$
where $l(A_M)$ is the minimal number of generators of $A_M$.
Moreover $(r,s)$ is the signature of $\bf U$,
and we assume that $r>p$ and $s>q$. By
Theorem 1.14.4 (analog of Witt's theorem) [Ni], it follows that there exists
a unique primitive embedding of $M$ into $\bf U$, up
to $O({\bf U})$. In particular, we have $f\circ i=i'\circ g$
for some $f\in O({\bf U})$ by uniqueness. Applying this to $M$, we
find that $fM=M'$, i.e. $M=M'~mod~O({\bf U})$.
This proves that (*) is injective. Likewise,
every $[L]\in\ol\cL_D(p,q)$ is represented by some $M\in\Pi$,
by the existence of primitive embedding $L\hra\bf U$.
This completes the proof. $\Box$

Note that the map (*) in the preceding proof is
well-defined without any of the hypotheses in the proposition.
However the example above, where $E=\R^{1,1}$
and $V\subset E$ a 1 dimensional positive definite,
shows that (*) need not be a bijection without those hypotheses.


\vfill\eject
\noi{\bf References:}\bs

\item{[Ca]} J.W.S. Cassels, {\it Rational quadratic forms}, Academic Press (1978).

\item{[CS]} J.H. Conway, and N.J.A. Sloane, {\it Sphere Packings,
Lattices and Groups}, 3rd Ed., Springer (1999).

\item{[GV]} S. Gukov, and C. Vafa,
{\it Rational Conformal Field Theories and Complex Multiplication},
hep-th/0203213.

\item{[HLOYI]} S. Hosono, B. Lian, K. Oguiso, and S.T. Yau,
{\it Counting Fourier-Mukai partners and applications},
math.AG/0202014.

\item{[HLOYII]} S. Hosono, B. Lian, K. Oguiso, and S.T. Yau,
{\it  Classification of c=2 Rational Conformal Field Theories
via the Gauss Product}, hep-th/0211230.

\item{[Mo]} G. Moore, {\it Attractors and Arithmetic}, hep-th/9807056;
 {\it Arithmetic and Attractors}, hep-th/9807087.

\item{[Ni]} V.V. Nikulin, {\it Integral symmetric bilinear forms and
some of their applications}, Math. USSR Izv. 14 (1980) 103-167.

\item{[Na]} K.S. Narain, {\it New Heterotic String Theories in Uncompactified Dimension $<10$},
Phys. Lett. B169 (1986) 41.

\item{[NSW]} K.S. Narain, M.H. Sarmadi, and E. Witten, {\it A note on toroidal compactification of heterotic
string theory}, Nucl. Phys. B279 (1987) 369.

\item{[Po]} J. Polchinski, {\it String Theory}, Vol. I,
Cambridge University Press (1998).

\item{[S]} J-P. Serre, {\it A Course in Arithmetic},
Springer (1973).

\item{[W]} K. Wendland
{\it Moduli Spaces of Unitary Conformal Field Theories},
Ph.D. thesis (available at
http://www-biblio.physik.uni-bonn.de/dissertationen/2000/doc/index.shtml).

\bs
\bs
\font\ninerm=cmr8 \font\sixrm=cmr6 \font\ninei=cmmi8 \font\sixi=cmmi6
\font\ninesy=cmsy8 \font\sixsy=cmsy6 \font\ninebf=cmbx8
\font\nineit=cmti8 \font\ninesl=cmsl8 \skewchar\ninei='177
\skewchar\sixi='177 \skewchar\ninesy='60 \skewchar\sixsy='60

\def\ninepoint{\def\rm{\fam0\ninerm}
\textfont0=\ninerm \scriptfont0=\sixrm \scriptscriptfont0=\fiverm
\textfont1=\ninei \scriptfont1=\sixi \scriptscriptfont1=\fivei
\textfont2=\ninesy \scriptfont2=\sixsy \scriptscriptfont2=\fivesy
\textfont\itfam=\ninei \def\it{\fam\itfam\nineit}\def\sl{\fam\slfam\ninesl}%
\textfont\bffam=\ninebf \def\bf{\fam\bffam\ninebf}\rm}
\def\footnotefont{\ninepoint}

{\footnotefont
\item{} Shinobu Hosono, Graduate School of Mathematical Sciences, University of Tokyo,
Komaba 3-8-1, Meguroku, Tokyo 153-8914, Japan. hosono@ms.u-tokyo.ac.jp

\item{} Bong H. Lian, Department of Mathematics, Brandeis University, Waltham, MA 02154, U.S.A.
lian@brandeis.edu

\item{} Keiji Oguiso, Graduate School of Mathematical Sciences, University of Tokyo,
Komaba 3-8-1, Meguroku, Tokyo 153-8914, Japan. oguiso@ms.u-tokyo.ac.jp

\item{} Shing-Tung Yau, Department of Mathematics, Harvard University,
Cambridge, MA 02138, U.S.A. yau@math.harvard.edu

\item{}
}
\end